\documentclass[12pt]{article}
\begin{document}

\setcounter{secnumdepth}{5}

\newcommand{\R}{I\hspace{-1.5 mm}R}
\newcommand{\N}{I\hspace{-1.5 mm}N}
\newcommand{\1}{I\hspace{-1.5 mm}I}
\newtheorem{proposition}{Proposition}[section]
\newtheorem{theorem}{Theorem}[section]
\newtheorem{lemma}[theorem]{Lemma}
\newtheorem{prop}[theorem]{Proposition}
\newtheorem{coro}[theorem]{Corollary}
\newtheorem{rem}[theorem]{Remark}
\newtheorem{ex}[theorem]{Example}
\newtheorem{claim}[theorem]{Claim}
\newtheorem{conj}[theorem]{Conjecture}
\newtheorem{remark}[theorem]{Remark}
\newtheorem{corollary}[theorem]{Corollary}


\title{Geometric inequalities via a general comparison principle for
interacting
  gases}
\author{M. Agueh\thanks{This paper was done while this author held a
postdoctoral
fellowship at the University of British Columbia.},\, N.
Ghoussoub\thanks{The
three
    authors were partially supported  by a grant from the Natural Science
and Engineering Research Council of Canada.}\,\, and X. Kang\thanks{This
paper
is part of this author's PhD's thesis under the supervision of N.
Ghoussoub.
}}
\date{Revised June 30, 2003}
\maketitle
\begin{center}
Pacific Institute for the Mathematical
Sciences\\and \\ Department of Mathematics,
The University of British Columbia\\
Vancouver, B. C.  V6T 1Z2, Canada
\end{center}

\section*{Abstract}

 The article builds on several recent advances in the  Monge-Kantorovich theory
of mass transport which have -- among other things -- led to new and quite
natural proofs for a wide range of geometric inequalities such as the ones
formulated by  Brunn-Minkowski, Sobolev, Gagliardo-Nirenberg, Beckner, Gross, 
Talagrand, Otto-Villani and their extensions by many others.  While this paper
continues in this spirit, we however propose here a basic framework to which
all of these  inequalities belong, and a general unifying principle from which
many of them follow. This basic inequality  relates  the relative total energy --
internal, potential and interactive -- of two arbitrary probability densities,
their Wasserstein distance, their barycentres and their
  entropy production
functional. The framework is remarkably encompassing as it implies many 
old geometric  -- Gaussian and Euclidean -- inequalities as well 
as new ones, while
allowing a direct and
unified way for computing best constants and extremals. As expected, such
inequalities also lead to exponential rates of convergence to equilibria
for solutions
of Fokker-Planck and McKean-Vlasov type equations.  The principle also leads to a
remarkable correspondence between ground state solutions of certain quasilinear
-- or semilinear -- equations  and stationary solutions of -- nonlinear -- 
Fokker-Planck type equations.

\tableofcontents

\section{Introduction}

The  recent advances in the  Monge-Kantorovich theory of mass transport  have --
among other things -- led to new and quite natural proofs for a wide range of
geometric inequalities. Most notable are  McCann's generalization of the 
Brunn-Minkowski's inequality \cite{mccann}, Barthe's proof of the -- reverse -- multidimensional Brascamp-Lieb inequality \cite{Ba}, Otto-Villani's
\cite{otto:generalization} and Cordero-Gangbo-Houdr\'e
\cite{cordero:inequalities} extensions of the Log Sobolev inequality
of Gross \cite{Gross} and Bakry-Emery \cite{bakry}, as well as 
Cordero-Nazaret-Villani's proof  \cite{cordero:mass}  of the Sobolev and the
Gagliardo-Nirenberg inequalities. We refer to the superb recent monograph
of Villani \cite{villani:topics} for more details on these remarkable
developments. 

This paper continues in this spirit, but our emphasis here is on developing
a framework for a  unified and compact approach to a substantial number of
these inequalities which originate in disparate areas of analysis and
geometry. 
The main idea is to try to describe the evolution of the total -- internal, potential
and interactive --  energy  of a system along an optimal transport that takes one
configuration to another,  taking
into account the  entropy production functional,  the
transport cost
(Wasserstein distance), as well as the displacement of their centres of mass. Once
this general
comparison principle is established, then  several -- new and old --
inequalities follow directly by simply considering different examples of --
admissible -- internal energies, and various 
confinement and interactive potentials. Others (e.g., Concentration of measure
phenomenon and Poincar\'e's inequality) will in turn follow from the well known
hierarchy between these inequalities.

Besides the obvious pedagogical relevance of a streamlining  approach,  we find it
interesting and intriguing that most of these inequalities appear as  different
manifestations of one basic principle in the theory of interacting gases that
compares the  energies  of two states of a system after one is transported ``at
minimal cost'' into another. Here is our framework which is already present in
McCann's thesis
\cite{mccann:thesis}.
 Let $\Omega$ be an open and convex subset of $\R^n$. The set of
probability densities over  $\Omega$ is
denoted by ${\cal P}_c(\Omega)=\{\rho: \Omega \to\R;\, \rho \geq 0 \,{\rm
and} \,
\int_\Omega\rho(x)dx=1\}$ and $\mbox{supp}\,\rho$ will stand for the support 
of $\rho\in{\cal P}_c(\Omega)$, that is the closure of $\{x\in\Omega:
\rho\neq 0\}$, while $|\Omega|$ will denote the Lebesgue measure of
$\Omega\subset\R^n$. 
Let $F:[0,\infty)\rightarrow \R$ be a differentiable function on
$(0,\infty)$, and let $V$ and $W$ be $C^2$-real valued functions on $\R^n$. 
The associated {\it Free Energy Functional} is
then defined on
${\cal P}_c(\Omega)$ as:
\[{\rm H}_V^{F,W}(\rho):=\int_{\Omega}\left[F(\rho) +\rho V
+\frac{1}{2}(W\star\rho)
\rho\right]\,\mbox{d}x,\]
which is the sum of the internal energy
${\rm H^F}(\rho):=\int_{\Omega}F(\rho) dx,$ 
the potential energy
$
{\rm H}_V(\rho):=\int_{\Omega} \rho V dx$
and the interaction energy
${\rm H}^{W}(\rho):=\frac{1}{2}\int_{\Omega}\rho(W\star\rho)\,\mbox{d}x.
$
Of importance is also the concept of {\it relative energy of $\rho_0$ with
respect to
$\rho_1$} simply defined as:
 ${\rm
H}^{F,W}_{V}(\rho_0|\rho_1):={\rm H}^{F,W}_{V}(\rho_0)- {\rm
H}^{F,W}_{V}(\rho_1),
$
where $\rho_0$ and $\rho_1$ are two probability densities.
  The {\it relative entropy production of $\rho$ with
respect to $\rho_V$} is normally defined as
\[
I_2(\rho|\rho_V)= \int_\Omega
\rho\Big|\,\nabla\left(F'(\rho)+V+W\star\rho)\right)\,\Big|^2\,\mbox{d}x
\]
in such a way that if $\rho_{_V}$ is a
probability density that satisfies
\[
       \nabla\left(F^\prime(\rho_{_V})+V+W\star \rho_{_V}\right) =
0\quad \mbox{a.e.}
\]
then
\[
I_2(\rho|\rho_{_V})= \int_\Omega \rho |\nabla\left(F^\prime
(\rho)-F'(\rho_{_V})+W\star (\rho-\rho_{_V}\right)|^2
\,\mbox{d}x.
\]
Our notation for the density $\rho_{_V}$ reflects this paper's emphasis on
its
dependence on the confinement potential, though it obviously also depends
on
$F$ and $W$.\\
We need the notion of Wasserstein distance $W_2$ between two probability
measures
$\rho_0$ and $\rho_1$ on $\R^n$, defined as:
\[
       W_2^2(\rho_0, \rho_1):=\inf_{\gamma \in \Gamma(\rho_0, \rho_1)}
\int_{\R^n\times\R^n}|x-y|^2 d\gamma (x,y),
\]
where $\Gamma(\rho_0, \rho_1)$ is the set of Borel probability measures
on $\R^n\times\R^n$ with marginals $\rho_0$ and
     $\rho_1$, respectively. 
The {\it barycentre} (or centre of mass) of a probability density $\rho$,
denoted
$
{\rm b}(\rho):=\int_{\R^n}x\rho(x)dx
$
will play a role in the presence of an interactive potential.

In this paper, we shall also deal with non-quadratic versions of the
entropy. For that we
call {\it Young function}, any strictly convex $C^1$-function
$c:\R^n\to\R$ such
that
$c(0)=0$ and $\lim_{|\,x\,|\rightarrow \infty}
\frac{c(x)}{|\,x\,|}=\infty$.  We denote by
$c^*$ its Legendre conjugate defined by
$
c^*(y)=\sup_{z\in \R^n}\{ y\cdot z-c(z)\}.
$
For any probability density $\rho$ on $\Omega$, we define the
{\it generalized relative entropy
production-type function of $\rho$ with respect to $\rho_V$ measured
against $c^*$} by
\[
{\cal I}_{c^*}(\rho|\rho_{_V}):=\int_\Omega \rho
c^\star\left(-\nabla\left(F^\prime(\rho)+V+W\star\rho\right)\right)\,\mbox{d}x,
\]
which is closely related to the {\it generalized relative entropy
production
function of  $\rho$ with respect to $\rho_V$  measured against $c^*$}
defined as:
\[
I_{c^*}(\rho|\rho_V):=\int_\Omega\rho\nabla\left(F'(\rho)+V+W\star\rho\right)
\cdot\nabla
c^\star\left(\nabla\left(F'(\rho)+V+W\star\rho\right)\right)\,\mbox{d}x.
\]
Indeed, the   convexity inequality
  $c^*(z)\leq z\cdot\nabla c^*(z)$ satisfied by any Young function
$c$,
  readily implies that
$
{\cal I}_{c^*}(\rho|\rho_V)\leq I_{c^*}(\rho|\rho_V).
$
Note that when $c(x)=\frac{|\,x\,|^2}{2}$, we have
\[ I_{c^*}(\rho|\rho_V)
=:I_2(\rho|\rho_V)=\int_\Omega
\rho\Big|\,\nabla\left(F'(\rho)+V+W\star\rho\right)\,\Big|^2\,\mbox{d}x=2{\cal
I}_{c^*}(\rho|\rho_V),\] and we denote ${\cal I}_{c^*}(\rho|\rho_V)$ by
${\cal
I}_2(\rho|\rho_V)$. 

Throughout this paper, the internal energy will be given
by a  differentiable function $F:[0,\infty)\rightarrow \R$  on
$(0,\infty)$ with $F(0)=0$ and $ x\mapsto x^nF(x^{-n})$ convex and
non-increasing. We denote by $P_F(x):=xF^\prime(x)-F(x)$   its
associated pressure function. The confinement potential will be given by 
a $C^2$-function $V:\R^n\to \R$ with $D^2V\geq \lambda I$, while the
interaction potential $W$ will be an even $C^2$-function  with $D^2W\geq \nu I$
where $\lambda,\nu \in \R$, and where $I$ stands for the identity map. 
 
In section \ref{sect2}, we start by establishing the following inequality relating the free
energies  of two  arbitrary probability densities, their Wasserstein distance,
their barycentres and their relative entropy production functional. The fact
that it yields many of the admittedly powerful geometric inequalities is
remarkable.\\  

\noindent{\bf Basic comparison principle for interactive gases:}
If $\Omega$ is any open, bounded and
convex subset of $\R^n$, then for any $\rho_0, \rho_1 \in {\cal P}_c(\Omega)$ 
satisfying
$\mbox{supp}\,\rho_0\subset \Omega$ and $P_F(\rho_0)\in
W^{1,\infty}(\Omega)$, and any Young function $c:\R^n\to\R$, we have:
\begin{equation}
\label{eqn:1}
  {\rm H}^{^{F,W}}_{_{V+c}}(\rho_0|\rho_1)+\frac{\lambda+\nu}{2}
W_2^2(\rho_0,
\rho_1)-\frac{\nu}{2}|{\rm b}(\rho_0)-{\rm b}(\rho_1)|^2
  \leq {\rm H}_{c+\nabla V\cdot x}^{^{-nP_F,2x\cdot\nabla W}}(\rho_0) +
{\cal 
I}_{c^*}(\rho_0 |\rho_V).
\end{equation}
Furthermore, equality holds in (\ref{eqn:1}) whenever
$\rho_0=\rho_1=\rho_{V+c}$,
where the latter satisfies
\begin{equation}
\label{eqn:2}
\nabla\left(F^\prime(\rho_{V+c}) +V+c +W\star\rho_{V+c}\right) =0 \quad
\mbox{a.e.}
\end{equation}
To give an idea about the strength of the above inequality, assume $V=W=0$ and
apply it with $\rho_0$ being any probability density $\rho$ satisfying 
$\mbox{supp}\,\rho\subset \Omega$ and $\rho_1=\rho_c$ the reference density. We
obtain:\\
  
\noindent {\bf The General Euclidean Sobolev Inequality:}
\begin{equation}
\label{eqn:7bis}
{\rm H}^{F+nP_F}(\rho) \leq
\int_\Omega\rho c^\star\left(-\nabla
(F^\prime\circ\rho)\right)\,\mbox{d}x +K_c,
\end{equation}
where $K_c$ is the unique constant determined by the equation
  \begin{equation}
\label{eqn:8bis}
F'(\rho_c)+c = K_c \, \, {\rm and} \, \,
\int_\Omega \rho_c=1.
\end{equation}
Applied to various -- displacement convex -- functionals $F$, we shall see in
section \ref{sect3} that (\ref{eqn:7bis}) already implies the Sobolev, the Gagliardo-Nirenberg and the
Euclidean $p$-Log Sobolev inequalities, allowing in the process a direct and
unified way for computing best constants and extremals. This formulation  also
points to an interesting fact: that the  various Sobolev inequalities are 
nothing but another manifestation of how free energy is controlled by entropy
production in appropriate systems.  

In section \ref{sect4}, we notice that inequality (\ref{eqn:1}) simplifies considerably
in the case where $c$ is a quadratic Young function of
the form $c(x):=c_\sigma(x)=\frac{1}{2\sigma}{|\,x\,|^2}$ for
$\sigma>0$, and we obtain: \\

\noindent{\bf The General Logarithmic Sobolev Inequality:}
 For all probability densities  $\rho_0$ and $\rho_1$ on $\Omega$,
satisfying
$\mbox{supp}\,\rho_0\subset \Omega$,  and $P_F(\rho_0)\in
W^{1,\infty}(\Omega)$,
we have for any $\sigma>0$,
  \begin{equation}
\label{eqn:3}
{\rm H}^{F,W}_{V}(\rho_0|\rho_1)+\frac{1}{2}(\lambda+\nu-\frac{1}{\sigma})
W_2^2(\rho_0,\rho_1) - \frac{\nu}{2}|{\rm b}(\rho_0)-{\rm
b}(\rho_1)|^2\leq
\frac{\sigma}{2}I_2(\rho_0|\rho_V).
\end{equation}
Minimizing the above inequality over $\sigma>0$ then yields:\\  

\noindent {\bf The HWBI
inequality
for interactive gases:}
\begin{equation}
\label{eqn:4}
{\rm H}^{F,W}_V(\rho_0|\rho_1)\leq
W_2(\rho_0,\rho_1)\sqrt{I_2(\rho_0|\rho_V)}
-\frac{\lambda+\nu}{2}W_2^2(\rho_0,\rho_1)+\frac{\nu}{2}|{\rm b}(\rho_0)-{\rm
b}(\rho_1)|^2.
\end{equation}
This extends the HWI inequality established in \cite{otto:generalization}
and
\cite{CMV}, with the additional ``B'' referring to the new barycentric
terms, and constitutes yet another extension of various powerful
inequalities by Gross \cite{Gross}, Bakry-Emery \cite{bakry}, Talagrand
\cite{T},
Otto-Villani \cite{otto:generalization}, Cordero \cite{cordero} and others. 

In section \ref{sect5}, we describe how these inequalities combined with the following
energy dissipation equation
\begin{equation}
\label{eqn:43}
\frac{d}{dt}\,{\rm H}^{F,W}_V\left(\rho(t)|\rho_V\right)
=-I_2\left(\rho(t)|\rho_V\right),
\end{equation}
provide  rates of convergence to
equilibria  for solutions to  
McKean-Vlasov type equations
\begin{equation}
\label{eqn:44}
\left\{\begin{array}{lcl}
\frac{\partial\rho}{\partial t} =
\mbox{div}\left\{\rho\nabla\left(F'(\rho)+V+W\star\rho\right)\right\}
&\mbox{in}&
(0,\infty)\times\R^n\\ \\
\rho(t=0)=\rho_0 &\mbox{in}& \{0\}\times\R^n.
\end{array}\right.
\end{equation}
One can then recover the recent results of Carrillo, McCann and Villani in
\cite{CMV}, which estimate the rate of convergence of solutions of (\ref{eqn:44}) to 
the equilibrium state.

In section \ref{sect6}, we apply inequality (\ref{eqn:1}) to the most basic system --
where no potential nor interaction energies are involved-- to obtain: \\

\noindent {\bf The Energy-Entropy production Duality Formula:} 
 For any probability density $\rho_0\in
W^{1,\infty}(\Omega)$ with support in $\Omega$,  and any $\rho_1\in {\cal
P}_a(\Omega)$, we have
\begin{equation}
\label{eqn:5}
-{\rm H}^F_{c}(\rho_1)\leq -{\rm H}^{F+nP_F}(\rho_0)+\int_\Omega \rho_0
c^\star\left(-\nabla(F^\prime\circ\rho_0)\right)\,\mbox{d}x.
\end{equation}
 Moreover, equality holds whenever $\rho_0=\rho_1=\rho_c$ where $\rho_c$ is
a probability
density on $\Omega$ such that $\nabla(F'(\rho_c)+c)=0$ a.e.\\
 Motivated by the recent work of
Cordero-Nazaret-Villani \cite{cordero:mass}, we show that
(\ref{eqn:5}) yields a statement of the following type:
\begin{equation}
\label{eqn:13}
\sup \{ J(\rho); \ \int_\Omega \rho(x)dx=1 \} \leq\inf \{
I(f);\  \int_\Omega
\psi(f(x))dx=1\},
\end{equation}
where
\begin{equation}
\label{eqn:14}
I(f)=\int_{\Omega} \left[c^*(-\nabla f (x)) -G\left(\psi\circ f(x)\right)\right] dx
\end{equation}
and
\begin{equation}
\label{eqn:15}
 J(\rho)= - \int_{\Omega} [ F(\rho(y)) + c(y) \rho(y) ] dy
\end{equation}
with $G(x)=(1-n)F(x)+nxF'(x)$ and where $\psi$ is computable from $F$ and
$c$.   Moreover, we have
equality in (\ref{eqn:13}) whenever there exists $\bar f$ (and
$\bar\rho=\psi(\bar
f)$)  that
satisfies the first order equation:
\begin{equation}
\label{eqn:16}
-(F' \circ \psi)'(\bar f) \nabla \bar f (x) = \nabla \mbox{$c(x)$ a.e.}
\end{equation}
In this case, the extrema are achieved at $\bar f$ (resp. $\bar \rho =
\psi (\bar f)$).
The latter is therefore a solution for the quasilinear (or semi-linear)
equation
     \begin{equation}
     \label{eqn:17}
\mbox{div} \{ \nabla c^{*} (-\nabla f) \} - (G\circ\psi)'(f)=\psi'
(f)
\end{equation}
since it is the $L^2$-Euler-Lagrange equation of $I$ on  
$
\{ f\in C^\infty_0(\Omega);\;\int_\Omega \psi(f(x))dx=1\}.
$
Equally interesting is the fact that $\psi (\bar f)$ is also
a stationary solution of the
(non-linear) Fokker-Planck equation:
\begin{equation}
    \label{eqn:18}
 \frac{\partial u}{\partial t}=   \mbox{div} \{ u \nabla
(F'(u)+c) \}
  \end{equation}
since $J$ is nothing but the Free Energy functional on  ${\cal
P}_a(\Omega)$, whose gradient
flow with respect  to the Wasserstein distance is precisely the evolution
equation (\ref{eqn:18}).
In other words, this is pointing to a remarkable correspondence between
Fokker-Planck  evolution
equations and certain quasilinear or semi-linear equations which appear as
Euler-Lagrange
equations of the entropy production functionals. 

 In conclusion to this introduction, we mention that this paper is an expanded
version of the unpublished but distributed manuscript \cite{agk1}. This
unifying and compact approach to so many important inequalities eventually led
us to make the paper as self-contained as possible so that it can serve as a
quick introduction to these basic tools of modern analysis. We should however warn
the reader that we have barely scratched the surface of the huge 
literature that exists on these basic inequalities, their various generalizations
and on the hierarchy and  relationships between them. Therefore, our
references are in no way complete nor exhaustive.  Fortunately many books and surveys
have already appeared on these topics and we refer the reader to the monograph of
Villani mentioned above, as well as to the book of Ledoux \cite{ledoux:topics} and
the recent survey of Gardner
\cite{gardner}.

 \section{Basic inequality between two configurations of interacting gases}
\label{sect2}
Here is our starting point.
\begin{theorem}
\label{theo2.1}
Let  $F:[0,\infty)\rightarrow \R$ be differentiable function on
$(0,\infty)$ with
$F(0)=0$ and $ x\mapsto x^nF(x^{-n})$
convex and non-increasing, and let $P_F(x):=xF^\prime(x)-F(x)$ be its
associated pressure
function. Let $V:\R^n\to \R$ be a $C^2$-confinement potential with
$D^2V\geq \lambda
I$,  and let $W$ be an even $C^2$-interaction potential with $D^2W\geq \nu
I$ where
$\lambda,\nu
\in \R$, and $I$ denotes the identity map. If $\Omega$ is any open, bounded
 and
convex subset of $\R^n$, then for any $\rho_0, \rho_1 \in {\cal P}_c(\Omega)$,
satisfying
$\mbox{supp}\,\rho_0\subset \Omega$ and $P_F(\rho_0)\in
W^{1,\infty}(\Omega)$, and any Young function $c:\R^n\to\R$, we have:
\begin{equation}
\label{eqn:2.1}
 {\rm H}^{^{F,W}}_{_{V+c}}(\rho_0|\rho_1)+\frac{\lambda+\nu}{2}
W_2^2(\rho_0,
\rho_1)-\frac{\nu}{2}|{\rm b}(\rho_0)-{\rm b}(\rho_1)|^2 
  \leq {\rm H}_{c+\nabla V\cdot x}^{^{-nP_F,2x\cdot\nabla W}}(\rho_0) +
+{\cal
I}_{c^*}(\rho_0 |\rho_V).
\end{equation}
Furthermore, equality holds in (\ref{eqn:2.1}) whenever
$\rho_0=\rho_1=\rho_{V+c}$,
where the latter satisfies
\begin{equation}
\label{eqn:2.2}
\nabla\left(F^\prime(\rho_{V+c}) +V+c +W\star\rho_{V+c}\right) =0 \quad
\mbox{a.e.}
\end{equation}
In particular, we have for any $\rho\in {\cal P}_c(\Omega)$ with
$\mbox{supp}\,\rho\subset\Omega$ and $P_F(\rho)\in W^{1,\infty}(\Omega)$,
{\footnotesize
\begin{equation}
\label{eqn:2.3}
 {\rm H}_{_{V-x\cdot \nabla V}}^{^{F+nP_F,\, W-2x\cdot\nabla
W}}(\rho)+
\frac{\lambda+\nu}{2} W_2^2(\rho,\rho_{V+c}) - \frac{\nu}{2}|{\rm
b}(\rho_0)-{\rm
b}(\rho_{V+c})|^2 \leq  
{\cal
I}_{c^*}(\rho |\rho_V)-{\rm
H}^{P_F,W}(\rho_{V+c})+K_{V+c},
 \end{equation}
}
where $K_{V+c}$ is a constant such that
\begin{equation}
\label{eqn:2.4}
F'(\rho_{V+c})+V+c+W\star \rho_{V+c}=K_{V+c} \, \, {\rm while} \,
\int_\Omega\rho_{V+c}=1.
\end{equation}

\end{theorem}
The proof is based on the recent advances in the theory of mass
transport as developed by Brenier \cite{brenier:polar}, Gangbo-McCann
\cite{gangbo:geometry}, \cite{GM}, Caffarelli \cite{caffareli:allocation}
and many others. For a survey, see Villani \cite{villani:topics}.  Here is
a brief summary of the needed results.\\
Fix a non-negative $C^1$, strictly convex function $d:\R^n\to \R$
such that $d(0)=0$. Given two probability measures $\mu$ and $\nu$ on
$\R^n$, the
minimum cost for transporting $\mu$ onto $\nu$ is given by
\begin{equation}
\label{eqn:2.5}
 W_d(\mu, \nu):=\inf_{\gamma \in \Gamma(\mu,
\nu)}\int_{\R^n\times\R^n}d(x-y)
d\gamma (x,y),
\end{equation}
where $\Gamma(\mu, \nu)$ is the set of Borel probability measures with
marginals
$\mu$ and $\nu$, respectively. When $d(x)=|\,x\,|^2$, we have that
$W_d=W^2_2$, where
$W_2$ is the Wasserstein distance. We say that a Borel map $T:\R^n\to
\R^n$ pushes $\mu$
forward to $\nu$,  if $\mu (T^{-1}(B))=\nu(B)$ for any Borel set
$B\subset\R^n$. The map
$T$ is then said to be $d$-optimal if
\begin{equation}
\label{eqn:2.6}
 W_d(\mu, \nu)
=\int_{\R^n}d(x-Tx)d\mu(x)=\inf_{S}\int_{\R^n}d(x-Sx)d\mu(x),
\end{equation}
where the infimum is taken over all Borel maps $S:\R^n\to \R^n$ that push
$\mu$
forward to $\nu$.
For quadratic cost functions $d(z)=\frac{1}{2}|z|^2$, Brenier
\cite{brenier:polar} characterized the optimal transport map $T$ as the
gradient
of a convex function. An analogous result holds for general cost functions
$d$, provided
convexity is replaced by an appropriate notion of $d$-concavity. See
\cite{gangbo:geometry}, \cite{caffareli:allocation} for details.

 Here is the lemma which leads to our main inequality (\ref{eqn:2.1}). It
is
essentially a compendium of various observations by several authors. It
describes the
evolution of a generalized energy functional along optimal transport. The
key idea  
behind it, is the concept of {\it displacement convexity} introduced by
McCann
\cite{mccann}. For generalized cost functions, and when $V=0$, it was
first obtained by
Otto \cite{otto:doubly} for the Tsallis entropy functionals and by Agueh
\cite{agueh:existence} in general. The case of a nonzero confinement
potential $V$ and an
interaction potential $W$ was included in \cite{cordero:inequalities},
\cite{CMV}. Here,
we state the results when the cost function is quadratic,
$d(x)=|\,x\,|^2$.

\begin{lemma}
\label{lem:1}
 Let $\Omega\subset\R^n$ be open, bounded and convex, and let $\rho_0$ and
$\rho_1$
be probability densities on $\Omega$, with
$\mbox{supp}\,\rho_0\subset\Omega$, and
$P_F(\rho_0)\in W^{1,\infty}(\Omega)$. Let $T$ be the  optimal map that
pushes $\rho_0\in{\cal
P}_a(\Omega)$ forward to $\rho_1\in{\cal P}_c(\Omega)$ for the quadratic
cost $d(x)=|\,x\,|^2$.
Then
\begin{trivlist}
\item 1) Assume $F:[0,\infty)\rightarrow \R$ is differentiable on
$(0,\infty)$,
$F(0)=0$ and $ x\mapsto x^nF(x^{-n})$ is convex and non-increasing, then
the internal energy satisfies:
\begin{equation}
\label{eqn:2.7}
{\rm H}^F(\rho_1)-{\rm H}^F(\rho_0) \geq \int_{\Omega}\rho_0(T-I)\cdot
\nabla\left(F'(\rho_0)\right)\,\mbox{d}x  .
\end{equation}
\item 2) Assume $V:\R^n\to \R$ is such that $D^2V\geq \lambda I$ for some
$\lambda\in \R$,
then the  potential energy satisfies
\begin{equation}
\label{eqn:2.8}
{\rm H}_V(\rho_1)-{\rm H}_V(\rho_0) \geq \int_{\Omega}\rho_0(T-I)\cdot
\nabla V\mbox{d}x+\frac{\lambda}{2}W^2_2(\rho_0,\rho_1).
\end{equation}
\item 3) Assume $W:\R^n\to \R$ is even, and $D^2W\geq \nu I$ for some $\nu\in
\R$, then
 the interaction energy satisfies
\begin{equation}
\label{eqn:2.9}
{\rm H}^W(\rho_1)-{\rm H}^W(\rho_0) \geq
\int_{\Omega}\rho_0(T-I)\cdot
\nabla (W\star \rho_0)\mbox{d}x    
  +\frac{\nu}{2}\left (W^2_2(\rho_0,\rho_1) -
|{\rm b}(\rho_0)-b(\rho_1)|^2\right).\nonumber
\end{equation}
\end{trivlist}
\end{lemma}
\noindent {\bf Proof:} If $T$ ($T=\nabla\psi$, where $\psi$ is convex) is the 
optimal map that pushes
$\rho_0\in{\cal
P}_a(\Omega)$ forward to $\rho_1\in{\cal P}_c(\Omega)$ for the quadratic
cost $d(x)=|\,x\,|^2$, one can then  define a path of probability
densities joining them, by letting $\rho_t$ be the push-forward measure of
$\rho_0$ by the map $T_t=(1-t)I + tT$. The key idea behind the estimate for
the internal energy is the fact first noticed by McCann
\cite{mccann}, that under the above assumptions on
$F$,  the function $t\mapsto {\rm H}^F(\rho_t)$ is convex on
$[0,1]$, which -- at least for smooth $\rho_t$
-- essentially leads to (\ref{eqn:2.7}) via the following
inequality for the internal energy:
\begin{equation}
\label{eqn:2.10}
{\rm H}^F(\rho_1)-{\rm H}^F(\rho_0) \geq\big[\frac{d}{dt} {\rm
H}^F(\rho_t)\big]_{t=0} = -\int_\Omega
F'(\rho_0)\,\mbox{div}\left(\rho_0(T-I)\right)\,\mbox{d}x.
\end{equation}
We shall use here another approach due to Agueh \cite{agueh:existence}
as it is more elementary and is applicable to other cost functions.\\ 
 First note that $\nabla T=\nabla \psi^2$ is diagonalizable with positive eigenvalues 
for $\rho_0$ a.e., and satisfies the
Monge-Amp\`ere equation 
\begin{equation}
\label{en1}
0\neq \rho_0(x)=\rho_1\left(T(x)\right)\mbox{det}\,\nabla T(x) \quad   \ \rho_0 \ 
\ {\rm a.e.}  
\end{equation}
So, $\rho_1\left(T(x)\right)\neq 0$ for $\rho_0$ a.e. 
Here, $\nabla T(x)=\nabla^2\psi(x)$ denotes the derivative in the sense of Aleksandrov of
$\psi$ (see McCann \cite{mccann}). 
Set 
$A(x)= x^n F(x^{-n})$, which is non-increasing by assumption, hence the pressure $P_F$ is
non-negative  and $x\mapsto\frac{F(x)}{x}$ is also non-increasing. Use that $F(0)=0$,
$T_{\#}\rho_0=\rho_1$ and (\ref{en1}), to obtain that
\begin{eqnarray}
\label{en2}
{\rm H}^F(\rho_1) = \int_{[\rho_1\neq 0]}\frac{F\left(\rho_1(y)\right)}
{\rho_1(y)}\rho_1(y)\,\mbox{d}y &=&
\int_{\Omega}\frac{F\left(\rho_1(Tx)\right)}{\rho_1(Tx)}\rho_0(x)\,\mbox{d}x \nonumber\\ &=&
\int_\Omega F\left(\frac{\rho_0(x)}{\mbox{det}\,\nabla T(x)}\right)\mbox{det}\,\nabla
T(x)\,\mbox{d}x.
\end{eqnarray}
Comparing the geometric mean $\left(\mbox{det}\,\nabla T(x)\right)^{1/n}$ with  
the arithmetic mean $\frac{\mbox{tr}\,\nabla T(x)}{n}$, we get 
$\frac{1}{\mbox{det}\,\nabla T(x)} \geq \left(\frac{n}{\mbox{tr}\,\nabla T(x)}\right)^n$,
and since $x\mapsto\frac{F(x)}{x}$ is non-decreasing, we obtain
\begin{equation}
\label{en3}
 F\left(\frac{\rho_0(x)}{det\,\nabla T(x)}\right)\,\mbox{det}\,\nabla T(x) \geq \Lambda^n F\left(\frac{\rho_0(x)}{\Lambda^n}\right) = \rho_0(x)A\left(\frac{\Lambda}{\rho_0(x)^{1/n}}\right),
\end{equation}
where $\Lambda := \frac{\mbox{tr}\,\nabla T(x)}{n}$. 
Next, we use that $A^\prime(x)=-nx^{n-1}P_F(x^{-n})$ and that $A$ is convex, to
obtain that
\begin{eqnarray}
\label{en4}
\rho_0(x)A\left(\frac{\Lambda}{\rho_0(x)^{1/n}}\right) &\geq& \rho_0(x)\left[\,A\left(\frac{1}{\rho_0(x)^{1/n}}\right)+A^\prime\left(\frac{1}{\rho_0(x)^{1/n}}\right)\left(\frac{\Lambda-1}{\rho_0(x)^{1/n}}\right)\right] \nonumber\\
&=& \rho_0(x)\left[\frac{F\left(\rho_0(x)\right)}{\rho_0(x)} - n(\Lambda-1)\frac{P_F\left(\rho_0(x)\right)}{\rho_0(x)}\right] \nonumber\\
&=& F\left(\rho_0(x)\right) - P_F\left(\rho_0(x)\right)\mbox{tr}\,(\nabla T(x)-I).
\end{eqnarray}
We combine (\ref{en2}) - (\ref{en4}), to conclude that
\begin{eqnarray}
{\rm H}^F(\rho_1)-{\rm H}^F(\rho_0) &\geq& -\int_\Omega P_F\left(\rho_0(x)\right)\,tr\,(\nabla T(x)-I)\,\mbox{d}x\nonumber\\
 &=& -\int_\Omega P_F\left(\rho_0(x)\right)\,\mbox{div}\,(T(x)-I)\,\mbox{d}x\nonumber\\
& \geq& \int_\Omega \rho_0\left(T-I\right)\cdot\nabla\left(F'(\rho_0)\right)\,\mbox{d}x.
\end{eqnarray} 

\noindent {\bf(2)} As noted in \cite{cordero:inequalities}, the fact that $D^2V\geq
\lambda I$, which means that
\[
V(b)-V(a) \geq \nabla V (a)\cdot(b-a) +\frac{\lambda}{2} |\,a-b\,|^2
\]
for all $a,b \in \R^n$,  easily implies (\ref{eqn:2.8}) via the following
inequality for the corresponding potential energy:

\begin{eqnarray}
\label{eqn:2.11}
{\rm H}_V(\rho_1)-{\rm H}_V(\rho_0) &\geq& \big[\frac{d}{dt} {\rm
H}_V(\rho_t)\big]_{t=0}+\frac{\lambda}{2}\int_{\Omega}|(T-I)(x)|^2\rho_0(x)dx\nonumber
\\
&=& -\int_\Omega V\mbox{div}\left(\rho_0(T-I)\right)\,\mbox{d}x +
\frac{\lambda}{2}W^2_2(\rho_0,\rho_1).
\end{eqnarray}

\noindent {\bf(3)} The following proof of (\ref{eqn:2.9}) appeared in 
Cordero-Gangbo-Houdr\'e \cite{cordero:inequalities}. Rewrite the interaction energy as follows:
\begin{eqnarray}
\label{eqn:2.12}
{\rm
H}^{W}(\rho_1)&=&\frac{1}{2}\int_{\Omega\times\Omega}W(x-y)\rho_1(x)\rho_1(y)\,\mbox{d}x\mbox{d}y\nonumber\\
&=&\frac{1}{2}\int_{\Omega\times\Omega}W(T(x)-T(y))\rho_0(x)\rho_0(y)\,\mbox{d}x\mbox{d}y\nonumber\\
&=&\frac{1}{2}\int_{\Omega\times\Omega}W\left(x-y+(T-I)(x)-(T-I)(y)\right)\rho_0(x)\rho_0(y)\,\mbox{d}x\mbox{d}y\nonumber\\
&\geq& \frac{1}{2}\int_{\Omega\times \Omega}\left[W(x-y)+\nabla
W(x-y)\cdot\left((T-I)(x)-(T-I)(y)\right)\rho_0(x)\rho_0(y)\right]\,\mbox{d}x\mbox{d}y\nonumber\\
&
&\;\;+\frac{\nu}{4}\int_{\Omega\times\Omega}|(T-I)(x)-(T-I)(y)|^2\rho_0(x)\rho_0(y)\,\mbox{d}x\mbox{d}y\nonumber\\
&=& {\rm H}^W(\rho_0)+\frac{1}{2}\int_{\Omega\times\Omega}\nabla
W(x-y)\cdot\left((T-I)(x)-(T-I)(y)\right)
\rho_0(x)\rho_0(y)\,\mbox{d}x\mbox{d}y\nonumber\\
& &\;\;+\frac{\nu}{4}\int_{\Omega\times\Omega}|(T-I)(x)-(T-I)(y)|^2
\rho_0(x)\rho_0(y)\,\mbox{d}x\mbox{d}y,
  \end{eqnarray}
where we used above that $D^2W\geq\nu I$. The last term of the subsequent
inequality can be written as:
\begin{eqnarray}
\label{eqn:2.13}
\lefteqn{\int_{\Omega\times\Omega}|(T-I)(x)-(T-I)(y)|^2\rho_0(x)\rho_0(y)\,\mbox{d}x\mbox{d}y}
\nonumber \\
& & = 2\int_{\Omega}|(T-I)(x)|^2\rho_0(x)\,\mbox{d}x
-2\Big|\int_{\R^n}(T-I)(x)\rho_0(x)\,\mbox{d}x\Big|^2\nonumber \\
& & = 2\int_{\Omega}|(T-I)(x)|^2\rho_0(x)\,\mbox{d}x -2|{\rm
b}(\rho_1)-{\rm b}(\rho_0)|^2.
\end{eqnarray}
 And since $\nabla W$ is odd (because $W$ is even), we get for the second
term of (\ref{eqn:2.12})
\begin{eqnarray}
\label{eqn:2.14}
\lefteqn{\int_{\Omega\times\Omega}\left[\nabla
W(x-y)\cdot\left((T-I)(x)-(T-I)(y)\right)\right]\rho_0(x)\rho_0(y)\,\mbox{d}x\mbox{d}y}\nonumber\\
& & = 2\int_{\Omega\times \Omega} \nabla W(x-y)\cdot
(T-I)(x)\rho_0(x)\rho_0(y)\,\mbox{d}x\mbox{d}y\nonumber\\
& &= 2\int_{\Omega\times \Omega} \rho_0 (T-I)\cdot \nabla (W\star \rho_0)
\,\mbox{d}x.
\end{eqnarray}
Combining (\ref{eqn:2.12}) - (\ref{eqn:2.14}), we obtain that
\begin{eqnarray*}
\lefteqn{{\rm H}^W(\rho_1)-{\rm H}^W(\rho_0)}\\
& &  \geq \int_{\Omega\times\Omega} \rho_0(T-I)\cdot \nabla (W\star
\rho_0)\,\mbox{d}x + \frac{\nu}{2}\left(\int_{\Omega}|(T-I)(x)|^2\rho_0dx
-|{\rm b}(\rho_0)-{\rm b}(\rho_1)|^2\right).
\end{eqnarray*}
This complete the proof of (\ref{eqn:2.9}).\\

\noindent{\bf Proof of Theorem \ref{theo2.1}:} Adding (\ref{eqn:2.7}),
(\ref{eqn:2.8}) and (\ref{eqn:2.9}), one gets
\begin{eqnarray}
\label{eqn:2.15}
\lefteqn{{\rm H}^{F,W}_V(\rho_0)- {\rm
H}^{F,W}_V(\rho_1)+\frac{\lambda+\nu}{2} W_2^2(\rho_0,\rho_1)
-\frac{\nu}{2}|{\rm b}(\rho_0)-{\rm b}(\rho_1)|^2} \\
& &\leq \int_\Omega
(x-Tx)\cdot\rho_0\nabla\left(F^\prime(\rho_0)+V+W\star\rho_0\right)\,\mbox{d}x.\nonumber
\end{eqnarray}
Since $\rho_0\nabla(F^\prime(\rho_0))=\nabla\left(P_F(\rho_0)\right)$,
we integrate by part
$\int_\Omega\rho_0\nabla\left(F^\prime(\rho_0)\right)\cdot x\,\mbox{d}x$,
and obtain that
\[\int_\Omega  x\cdot\nabla(F^\prime(\rho_0)+V+W\star \rho_0)\rho_0={\rm
H}_{x\cdot\nabla V}^{^{-nP_F,\,2x\cdot \nabla W}}(\rho_0).\]
This leads to
\begin{eqnarray}
\label{eqn:2.16}
\lefteqn{{\rm H}^{^{F,W}}_V(\rho_0)-{\rm H}^{F,W}_V(\rho_1)
+\frac{\lambda+\nu}{2} W_2^2(\rho_0, \rho_1) -\frac{\nu}{2}|{\rm
b}(\rho_0)-{\rm b}(\rho_1)|^2} \\
& &\leq \quad {\rm H}_{x\cdot\nabla V}^{^{-nP_F, \, 2x\cdot \nabla
W}}(\rho_0)-\int_\Omega\rho_0\nabla\left(F^\prime(\rho_0)+V+W\star\rho_0\right)\cdot
T(x)\, \mbox{d}x.\nonumber
\end{eqnarray}
Now, use Young's inequality  to get
\begin{eqnarray}
\label{eqn:2.17}
\lefteqn{-\nabla\left(F^\prime\left(\rho_0(x)\right) +
V(x)+(W\star\rho_0)(x)\right)\cdot T(x)}\\
& & \leq c\left(T(x)\right) +
c^\star\left(-\nabla\left(F^\prime(\rho_0(x))+V(x)+(W\star\rho_0)(x)
\right)\right),\nonumber
\end{eqnarray}
and deduce that 
\begin{eqnarray}
\label{eqn:2.18}
\lefteqn{{\rm H}^{F,W}_V(\rho_0)-{\rm
H}^{F,W}_V(\rho_1)+\frac{\lambda+\mu}{2} W_2^2(\rho_0,
\rho_1)-\frac{\nu}{2}|{\rm b}(\rho_0)-{\rm b}(\rho_1)|^2}\\
& &\leq {\rm H}_{x\cdot\nabla V}^{-nP_F, 2x\cdot \nabla W}(\rho_0) +
\int_\Omega \rho_0
c^\star\left(-\nabla\left(F^\prime(\rho_0)+V+W\star\rho_0)\right)\right)+\int_\Omega
c(Tx) \rho_0\,\mbox{d}x.\nonumber
\end{eqnarray}
\noindent Finally, use again that $T$ pushes $\rho_0$ forward to $\rho_1$,
to rewrite
 the last integral on the right hand side of (\ref{eqn:2.18}) as
$\int_\Omega c(y)
\rho_1(y) dy$  to obtain (\ref{eqn:2.1}).\\ Now, set
$\rho_0=\rho_1:=\rho_{V+c}$ in
(\ref{eqn:2.16}). We have that $T=I$, and equality then holds in
(\ref{eqn:2.16}).
Therefore, equality holds in (\ref{eqn:2.1}) whenever equality holds in
(\ref{eqn:2.17}),
where $T(x)=x$. This occurs when (\ref{eqn:2.2}) is satisfied.\\
\noindent(\ref{eqn:2.3}) is straightforward when choosing $\rho_0:=\rho$
and
$\rho_1:=\rho_{V+c}$ in (\ref{eqn:2.1}). 

\section{The General Euclidean Sobolev Inequality}
\label{sect3}

We start with the following general inequality, which can be seen as an 
extension of the various Euclidean Sobolev inequalities, since once applied to
appropriate functionals $F$ and $c$, one gets the Sobolev, the
Gagliardo-Nirenberg and the Euclidean $p$-Log Sobolev inequalities. 

\begin{theorem} {\rm (The General Sobolev Inequality)}
\label{theo3.1}
Under the hypothesis of Theorem 2.1, assume that $V$ and $W$ are also convex. 
Then,  for any {\sl Young} function
$c:\R^n\to\R$, and any  $\rho\in {\cal P}_c(\Omega)$ with
$\mbox{supp}\,\rho\subset\Omega$ and $P_F(\rho)\in W^{1,\infty}(\Omega)$, we have 
\begin{equation}
\label{eqn:2.3.1}
 {\rm H}_{_{-V^*(\nabla V)}}^{^{F+nP_F,\, W-2x\cdot\nabla
W}}(\rho) \nonumber \leq \int_\Omega\rho
c^\star\left(-\nabla\left(F'(\rho)+V+W\star\rho\right)\right)\,\mbox{d}x-{\rm
H}^{P_F,W}(\rho_{V+c})+K_{V+c},
\end{equation}
where $\rho_{V+c}$ is the probability density and  $K_{V+c}$ is the
constant satisfying
 \begin{equation}
\label{eqn:2.4.1}
F'(\rho_{V+c})+V+c+W\star \rho_{V+c}=K_{V+c}.
\end{equation}
In particular, if $V=W=0$, we have 
\begin{equation}
\label{eqn:7}
{\rm H}^{F+nP_F}(\rho) \leq
\int_\Omega\rho c^\star\left(-\nabla
(F^\prime\circ\rho)\right)\,\mbox{d}x +K_c,
\end{equation}
where $K_c$ is the unique constant determined by the equation
  \begin{equation}
\label{eqn:8}
F'(\rho_c)+c = K_c \, \, {\rm and} \, \,
\int_\Omega \rho_c=1.
\end{equation}
\end{theorem}
\noindent{\bf Proof:} This follows immediately from inequality (\ref{eqn:2.3})
in Theorem 2.1. Indeed, if $\lambda+\nu\geq 0$, then the term involving the
Wasserstein distance can be omitted from the equation, while if $W$ is convex,
then the
barycentric term can also be omitted. If $V$ is strictly convex, then $V-x\cdot
\nabla V=-V^*(\nabla V)$. \\ 
Now if $V=W=0$, we obtain the remarkably simple inequality:
\begin{equation}
\label{eqn:9}
{\rm H}^{F+nP_F}(\rho) \leq
\int_\Omega\rho c^\star\left(-\nabla
(F^\prime\circ\rho)\right)\,\mbox{d}x-{\rm H}^{P_F}(\rho_c)+K_c,
\end{equation}
where $K_c$ is the unique constant determined by (\ref{eqn:8}). 
Finally, we obtain (\ref{eqn:7})  by noting that ${\rm H}^{P_F}(\rho_c)$ is
always non-negative. 

\subsection{Euclidean Log-Sobolev inequalities}
 
The following optimal Euclidean $p$-Log Sobolev inequality was first
established by Beckner in \cite{beckner:geometry} for $p=1$, and  by
Del-Pino and Dolbeault \cite{deldol:optimal} for $1< p < n$. The case 
where $p>n$ was
established
recently and independently by I. Gentil \cite {gentil} who used the
Pr\'ekopa-Leindler inequality
and the Hopf-Lax semi-group associated to the Hamilton-Jacobi equation.

 \begin{corollary} ({\rm General Euclidean Log-Sobolev inequality}) 
\label{coro4.7}
 Let $\Omega\subset\R^n$ be open bounded and convex, and let
$c:\R^n\rightarrow \R$ be a Young functional such that its conjugate
$c^\star$ is $p$-homogeneous for some $p>1$. Then,
\begin{equation}
\label{eqn:3.38}
\int_{\R^n} \rho\ln\rho\,\mbox{d}x \leq \frac{n}{p}\ln\left(\frac{p}
{ne^{p-1}\sigma_c^{p/n}}\int_{\R^n} \rho c^\star\left(-\frac{\nabla
\rho}{\rho}\right)\,\mbox{d}x\right),
\end{equation}
for all probability densities $\rho$ on $\R^n$, such that
$\mbox{supp}\,\rho\subset\Omega$ and $\rho\in W^{1,\infty}(\R^n)$. Here,
$\sigma_c:=\int_{\R^n} e^{-c}\,\mbox{d}x$.  Moreover, equality holds in
(\ref{eqn:3.38}) if $\rho(x)=K_\lambda e^{-\lambda^{q}c(x)}$ for some
$\lambda>0$,  where
$K_\lambda=\left(\int_{\R^n}e^{-\lambda^qc(x)}\,\mbox{d}x\right)^{-1}$
and $q$ is the conjugate of $p$ ($\frac{1}{p}+\frac{1}{q}=1)$.

\end{corollary}
{\bf Proof:} Use $F(x)=x\ln(x)$ and $V=W=0$   in
(\ref{eqn:2.3}). Note that $P_F(x)=x$, and then, $H^{P_F}(\rho)=1$ for any
$\rho\in {\cal P}_c(\R^n)$. So,
$\rho_c(x)=\frac{e^{-c(x)}}{\sigma_c}$. We then have for
     $\rho\in{\cal P}_c(\R^n)\cap W^{1,\infty}(\R^n)$ such that
$\,\mbox{supp}\,\rho\subset\Omega$,
\begin{equation}
\label{eqn:3.39}
\int_\Omega \rho\ln \rho \,\mbox{d}x  \leq \int_{\R^n}
\rho c^\star\left(-\frac{\nabla\rho}{\rho}\right)\,\mbox{d}x
-n-\ln\left(\int_{\R^n}e^{-c(x)}\,\mbox{d}x\right),
\end{equation}
with equality when $\rho=\rho_c$.

Now assume that $c^\star$ is $p$-homogeneous and set  $\Gamma^c_\rho
=\int_{\R^n}
\rho c^\star\left(-\frac{\nabla \rho}{\rho}\right)\,\mbox{d}x$. Using
$c_\lambda(x):=c(\lambda x)$ in (\ref{eqn:3.39}), we get for $\lambda>0$
that
\begin{equation}
\label{eqn:3.40}
\int_{\R^n} \rho\ln \rho \,\mbox{d}x \leq \int_{\R^n} \rho c^\star
\left(-\frac{\nabla\rho}{\lambda\rho}\right)\,\mbox{d}x+n\ln\lambda-n-\ln\sigma_c,
\end{equation}
for all $\rho\in{\cal P}_c(\R^n)$ satisfying
$\mbox{supp}\,\rho\subset\Omega$
 and $ \rho\in W^{1,\infty}(\Omega)$. Equality holds in (\ref{eqn:3.40})
if $\rho_\lambda
(x)=\left(\int_{\R^n}e^{-\lambda^q
c(x)}\,\mbox{d}x\right)^{-1}e^{-\lambda^q c(x)}$. Hence
\[\int_{\R^n} \rho\ln\rho \,\mbox{d}x \leq -n-\ln
\sigma_c+\inf_{\lambda>0}
\left(G_\rho(\lambda)\right),\]
      where
\[G_\rho(\lambda)=n\ln(\lambda) + \frac{1}{\lambda^p}\int_{\R^n} \rho
c^\star\left(-\frac{\nabla\rho}{\rho}\right) = n\ln(\lambda) +
\frac{\Gamma^c_{\rho}}{\lambda^p}.\] The infimum of $G_\rho(\lambda)$
over $\lambda>0$ is
attained at
$\bar{\lambda}_\rho=\left(\frac{p}{n}\Gamma^c_{\rho}\right)^{1/p}$.
Hence
\begin{eqnarray*}
\int_{\R^n}  \rho\ln\rho \,\mbox{d}x &\leq&
G_\rho(\bar{\lambda}_\rho)-n-\ln(\sigma_c)\\
&=& \frac{n}{p}\ln\left(\frac{p}{n}\Gamma^c_{\rho}\right)
+\frac{n}{p}-n-\ln(\sigma_c)\\
&=& \frac{n}{p}\ln\left(\frac{p}{n
e^{p-1}\sigma_c^{p/n}}\Gamma^c_{\rho}\right),
\end{eqnarray*}
for all probability densities $\rho$ on $\R^n$, such that $\mbox{supp}\,
\rho\subset\Omega$, and $\rho\in W^{1,\infty}(\R^n)$.

\begin{corollary} ({\rm Optimal Euclidean $p$-Log Sobolev inequality})\\
   \label{coro4.8}
\begin{equation}
\label{eqn:3.41}
\int_{\R^n} |\,f\,|^p\ln(|\,f\,|^p)\,\mbox{d}x \leq
\frac{n}{p}\ln\left(C_p\int_{\R^n} |\,\nabla f\,|^p\,\mbox{d}x\right),
\end{equation}
holds for all $p\geq 1$, and for all $f\in W^{1,p}(\R^n)$ such that
$\|\,f\,\|_p=1$, where
\begin{equation}
\label{eqn:3.42}
C_p:=\left\{\begin{array}{lcl}
\left(\frac{p}{n}\right)\left(\frac{p-1}{e}\right)^{p-1}
\pi^{-\frac{p}{2}}\left[\frac{\Gamma(\frac{n}{2}+1)}{\Gamma(\frac{n}{q}+1)}
\right]^{\frac{p}{n}}
&\mbox{if}& p>1,\\ \\
\frac{1}{n\sqrt{\pi}}
\left[\Gamma(\frac{n}{2}+1)\right]^{\frac{1}{n}} &\mbox{if}&  p=1,
\end{array}\right.
\end{equation}
and $q$ is the conjugate of $p$ $(\frac{1}{p}+\frac{1}{q}=1)$.\\
For $p>1$, equality holds in (\ref{eqn:3.41}) for
$f(x)=Ke^{-\lambda^{q}\frac{|\,x-\bar{x}\,|^q}{q}}$ for some $\lambda>0$
and $\bar{x}\in\R^n$, where $K=\left(\int_{\R^n}e^{-(p-1)|\,\lambda
x\,|^q}\,\mbox{d}x\right)^{-1/p}$.
\end{corollary}

\noindent{\bf Proof:} First assume that $p>1$, and  set
$c(x)=(p-1)|\,x\,|^q$  and $\rho=|\,f\,|^p$ in (\ref{eqn:3.38}), where
$f\in C^\infty_c(\R^n)$ and
$\|\,f\,\|_p=1$. We have that $c^\star(x)=\frac{|\,x\,|^p}{p^p}$, and
then, $\int_{\R^n}\rho c^*\left(-\frac{\nabla\rho}{\rho}\right)\,\mbox{d}x
=\int_{\R^n}|\,\nabla f\,|^p\,\mbox{d}x$. Therefore, (\ref{eqn:3.38})
reads as
\begin{equation}
\label{eqn:3.43}
\int_{\R^n}|\,f\,|^p\ln(|\,f\,|^p)\,\mbox{d}x \leq
\frac{n}{p}\ln\left(\frac{p}{n e^{p-1}\sigma_c^{p/n}}\int_{\R^n}|\,\nabla
f\,|^p\,\mbox{d}x\right).
\end{equation}
Now, it suffices to note that
\begin{equation}
\label{eqn:3.44}
\sigma_c:=\int_{\R^n}e^{-(p-1)|\,x\,|^q}\,\mbox{d}x=\frac{\pi^{\frac{n}{2}}
\Gamma\left(\frac{n}{q}+1\right)}{(p-1)^{\frac{n}{q}}\Gamma\left(\frac{n}{2}+1\right)}.
\end{equation}

To prove the case where $p=1$, it is sufficient to apply the above to
$p_\epsilon=1+\epsilon$ for some arbitrary $\epsilon>0$. Note that
\[C_{p\epsilon}=\left(\frac{1+\epsilon}{n}\right)\left(\frac{\epsilon}{e}
\right)^{\epsilon}
\pi^{-\frac{1+\epsilon}{2}}\left[\frac{\Gamma
(\frac{n}{2}+1)}{\Gamma(\frac{n\epsilon}{1+\epsilon}+1)}\right]^{\frac{1+\epsilon}{n}},\]
so that when $\epsilon$ go to $0$, we have
\[\lim_{\epsilon\rightarrow 0}
C_{p_\epsilon}=\frac{1}{n\sqrt{\pi}}
\left[\Gamma\left(\frac{n}{2}+1\right)\right]^{\frac{1}{n}}=C_1.\]

\subsection{Sobolev and Gagliardo-Nirenberg inequalities}

\begin{corollary} ({\rm Gagliardo-Nirenberg inequalities}) 
\label{coro4.9}
Let $1<p<n$ and $r\in\left(0,\frac{np}{n-p}\right)$ such that $r\neq p$.
Set $\gamma:=\frac{1}{r}+\frac{1}{q}$, where
$\frac{1}{p}+\frac{1}{q}=1$. Then, for any $f \in W^{1,p}(\R^n)$  we have
\begin{equation}
\label{eqn:3.46}
\|f\|_r\leq C(p,r)\|\nabla f\|_p^\theta\,\|f\|_{r\gamma}^{1-\theta},
  \end{equation}
where $\theta$ is given by
\begin{equation}
\label{eqn:3.49}
\frac{1}{r} =\frac{\theta}{p^*}+\frac{1-\theta}{r\gamma},
\end{equation}
$p^* = \frac{np}{n-p}$ and where the best constant $C(p,r)>0$ can be
obtained by scaling.
\end{corollary} 

\noindent{\bf Proof:} Let $F(x)=\frac{x^\gamma}{\gamma-1}$,
where $1\neq \gamma > 1-\frac{1}{n}$, which follows from the fact that
$p\neq r\in
\left(0,\frac{np}{n-p}\right)$. For this value of $\gamma$, the function
$F$  satisfies the
conditions of Theorem \ref{theo3.1}.  Let $c(x)=\frac{r\gamma}{q}|\,
x\,|^q$ so that
$c^*(x)=\frac{1}{p(r\gamma)^{p-1}}|\, x \,|^p$, and set $V=W=0$.
Inequality (\ref{eqn:2.3})
then gives for all $f\in C^\infty_c(\R^n)$ such that $\|\,f\,\|_r=1$,
\begin{equation}
\label{eqn:3.50}
\left(\frac{1}{\gamma-1}+n\right) \int_{\R^n} |\,f\,|^{r\gamma} \leq
\frac{r\gamma}{p} \int_{\R^n} |\,\nabla f\,|^p
-H^{P_F}(\rho_\infty)+C_\infty.
\end{equation}
where $\rho_\infty=h_\infty^r$ satisfies
\begin{equation}
\label{eqn:3.51}
-\nabla h_\infty(x)= x|\,x\,|^{q-2}h^{\frac{r}{p}}(x)\;\;\mbox{a.e.},
\end{equation}
and where $C_\infty$ insures that $\int h_\infty^r=1$. The constants on
the
  right hand
side of (\ref{eqn:3.50}) are not easy to calculate, so one can obtain
$\theta$ and the best constant by a standard scaling procedure. Namely,
write (\ref{eqn:3.50}) as
\begin{equation}
\label{eqn:3.52}
\frac{r\gamma}{p}\frac{\|\nabla
f\|_{p}^{p}}{\|f\|_{_r}^{^p}}-\left(\frac{1}{\gamma-1}+n\right)\frac{\|f\|_{r\gamma}^{r\gamma}}{{\|f\|_{_r}}}\geq
H^{P_F}(\rho_\infty)-C_\infty=: C,
\end{equation}
for some constant $C$. Then apply (\ref{eqn:3.52}) to
$f_\lambda(x)=f(\lambda x)$ for $\lambda>0$. A minimization over $\lambda$
gives the required constant.\\

The limiting case where $r$ is the critical Sobolev exponent
$r=p^*=\frac{np}{n-p}$ (and then $\gamma=1-\frac{1}{n}$) leads to the
Sobolev
inequalities:

\begin{coro} ({\rm Sobolev inequalities}) 
\label{coro4.10}
If $1<p<n$, then for any $f \in W^{1,p}(\R^n)$,
\begin{equation}
\label{eqn:3.53}
  \|\,f\,\|_{p^*}\leq C(p,n)\|\,\nabla f\,\|_p
  \end{equation}
for some constant $C(p,n)>0$.
\end{coro}
\noindent{\bf Proof:} It follows directly from (\ref{eqn:3.50}), by using
$\gamma=1-\frac{1}{n}$ and $r=p^*$. 
Note that the scaling argument cannot be used here to compute the best
constant $C(p,n)$ in (\ref{eqn:3.53}), since $\|\,\nabla
f_\lambda\,\|_p^p=\lambda^{p-n}\|\,\nabla f\,\|^p_p$ and
$\|\,f_\lambda\,\|_r^p=\lambda^{p-n}\|\,f\,\|^p_r$ scale the same way in
(\ref{eqn:3.52}). Instead, one can proceed directly from (\ref{eqn:3.50})
to have that
\[\|\,f\,\|_{p^*}=1\leq\left(\frac{r\gamma}{p\left[H^{P_F}(\rho_\infty)-
C_\infty\right]}\right)^{1/p}\|\,\nabla f\,\|_p =
\left(\frac{p^*(n-1)}{np\left[H^{P_F}(\rho_\infty)-
C_\infty\right]}\right)^{1/p}\|\,\nabla
f\,\|_p ,\] which shows that
\begin{equation}
\label{eqn:3.54}
C(p,n)=\left(\frac{p^*(n-1)}{np\left[H^{P_F}(\rho_\infty)-C_\infty\right]}
\right)^{1/p},
\end{equation}
where $\rho_\infty=h_\infty^{p^*}=\left(\frac{p^*}{nq}|\,x\,|^q-
\frac{C_\infty}{n-1}\right)^{-n}$ is obtained from (\ref{eqn:3.51}), and
$C_\infty$
can be found using that $\rho_\infty$ is a probability density,
\begin{equation}
\label{eqn:3.55}
C_\infty=(1-n)\left[\int_{\R^n}\left(\frac{p^*}{nq}
|\,x\,|^q+1\right)^{-n}\,\mbox{d}x\right]^{p/n}.
\end{equation}

\section{The General Logarithmic Sobolev Inequality}
\label{sect4}

In this section, we consider the case where $c$ is a quadratic Young function of
the form $c(x):=c_\sigma(x)=\frac{1}{2\sigma}{|\,x\,|^2}$ for
$\sigma>0$. In this case, our basic inequality (1) simplifies considerably to
yield Theorem 4.1 below, which relates  the total energy of two arbitrary
probability densities, their Wasserstein distance, their barycentres and their
entropy production functional. This gives yet another remarkable  extension of
various powerful inequalities by Gross \cite{Gross}, Bakry-Emery\cite{bakry},
Talagrand \cite{T}, Otto-Villani \cite{otto:generalization}, 
Cordero\cite{cordero} and others.
\begin{theorem} {\rm (General Logarithmic Sobolev Inequality)}
\label{prop4.1}
Under the hypothesis of Theorem 2.1, 
we have for all $\rho_0, \rho_1\in {\cal P}_c(\Omega)$,
satisfying
$\mbox{supp}\,\rho_0\subset \Omega$,  and $P_F(\rho_0)\in
W^{1,\infty}(\Omega)$, and any
$\sigma>0$,
\begin{equation}
\label{eqn:3.16}
{\rm H}^{F,W}_{U}(\rho_0|\rho_1)+\frac{1}{2}(\mu+\nu-\frac{1}{\sigma})
W_2^2(\rho_0,\rho_1) - \frac{\nu}{2}|{\rm b}(\rho_0)-{\rm
b}(\rho_1)|^2\leq
\frac{\sigma}{2}I_2(\rho_0|\rho_U).
\end{equation}
 \end{theorem}
\noindent {\bf Proof:} Apply inequality (\ref{eqn:2.1}) with a quadratic Young
functional 
$c(x)=\frac{1}{2\sigma}|\,x\,|^2 $,
$V=U-c$ and $\lambda=\mu-\frac{1}{\sigma}$ to obtain
\begin{eqnarray}
\label{eqn:3.17}
{\rm H}_U^{F,W}(\rho_0|\rho_1)&+&\frac{1}{2}(\mu+\nu-\frac{1}{\sigma})
W_2^2(\rho_0,\rho_1)-\frac{\nu}{2}|{\rm b}(\rho_0)-{\rm b}(\rho_1)|^2 \\
  &\leq & {\rm H}^{-nP_F, 2x\cdot \nabla W}_{c+\nabla(U-c)\cdot
x}(\rho_0)+
\int_\Omega \rho_0 c^*\left(-\nabla\left(F'(\rho_0)+U-c+W\star
\rho_0\right)\right)
\,\mbox{d}x.\nonumber
\end{eqnarray}
Now we show the identity:
\[
{\cal I}_{c_\sigma^*}(\rho_0|\rho_V)+H^{-nP_F,2x\cdot\nabla 
W}_{c_\sigma+x\cdot\nabla
V}(\rho_0)={\cal
I}_{c^*_\sigma}(\rho_0|\rho_{V+c_\sigma})=
\frac{\sigma}{2}I_2(\rho_0|\rho_{V+c_\sigma}). 
\]
Indeed, by elementary computations, we have
\begin{eqnarray*}
\lefteqn{\int_\Omega \rho_0 c^*\left(-\nabla\left(F'\circ\rho_0+U-c+
W\star\rho_0\right)\right)\,\mbox{d}x}\nonumber\\
& & = \frac{\sigma}{2}\int_\Omega\rho_0\Big|\,\nabla\left(F'(\rho_0)+
U+W\star \rho_0\right)\,\Big|^2\,\mbox{d}x +
\frac{1}{2\sigma}\int_\Omega\rho_0|\,x\,|^2\,\mbox{d}x -\int_\Omega\rho_0
x\cdot\nabla\left(F'(\rho_0)\right)\,\mbox{d}x\nonumber\\ & & \quad
-\int_\Omega\rho_0
x\cdot\nabla U\,\mbox{d}x-\int_\Omega\rho_0 x\cdot\nabla (W\star
\rho_0)\,\mbox{d}x,\nonumber
\end{eqnarray*}
and
\[
{\rm H}^{-nP_F,2x\cdot\nabla W}_{c+\nabla(U-c)\cdot x}(\rho_0) = -{\rm
H}^{nP_F}(\rho_0)+\int_\Omega\rho_0 x\cdot\nabla
(W\star\rho_0)\,\mbox{d}x+\int_\Omega\rho_0 x\cdot\nabla
U\,\mbox{d}x-\frac{1}{2\sigma}\int_\Omega|\,x\,|^2\rho_0\,\mbox{d}x.
\]
By combining the last 2 identities, we can rewrite the right hand side
of (\ref{eqn:3.17}) as
\begin{eqnarray}
\label{eqn:3.18}
\lefteqn{{\rm H}^{-nP_F,2x\cdot\nabla W}_{c+\nabla\left(U-c\right)\cdot
x}(\rho_0)+\int_\Omega \rho_0 c^*\left(-\nabla(F'\circ\rho_0+U-c+W\star
\rho_0)\right)\,\mbox{d}x}\nonumber\\ & & =
\frac{\sigma}{2}\int_\Omega\rho_0|\,\nabla\left(F'(\rho_0)+U+W\star
\rho_0\right)\,|^2\,\mbox{d}x - \int_\Omega \rho_0
x\cdot\nabla\left(F'\circ\rho_0\right)\,\mbox{d}x-\int_\Omega
nP_F(\rho_0)\,\mbox{d}x\nonumber\\
 & &=
\frac{\sigma}{2}\int_\Omega\rho_0|\,\nabla\left(F'(\rho_0)+U+W\star\rho_0\right)\,|^2\mbox{d}x+\int_\Omega\mbox{div}\,(\rho_0
x)F'(\rho_0)\,\mbox{d}x-\int_\Omega nP_F(\rho_0)\,\mbox{d}x\nonumber\\
& &= \frac{\sigma}{2}\int_\Omega
\rho_0\Big|\,\nabla\left(F'(\rho_0)+U+W\star
\rho_0\right)\,\Big|^2\,\mbox{d}x + n\int_\Omega\rho_0
F'(\rho_0)\,\mbox{d}x +\int_\Omega x\cdot\nabla F(\rho_0)\,\mbox{d}x
\nonumber\\
& & \quad -\int_\Omega nP_F(\rho_0)\,\mbox{d}x\nonumber\\
& &
=\frac{\sigma}{2}\int_\Omega\rho_0\Big|\,\nabla\left(F'(\rho_0)+U+W\star
\rho_0\right)\,\Big|^2\mbox{d}x+\int_\Omega x\cdot\nabla
F(\rho_0)\,\mbox{d}x + n\int_\Omega F\circ\rho_0\,\mbox{d}x\nonumber\\
& &
=\frac{\sigma}{2}\int_\Omega\rho_0\Big|\,\nabla\left(F'(\rho_0)+U+W\star
\rho_0\right)\,\Big|^2\,\mbox{d}x.
\end{eqnarray}
Inserting (\ref{eqn:3.18}) into (\ref{eqn:3.17}), we conclude
(\ref{eqn:3.16}). 

\subsection{HWBI inequalities}
We now establish the HWBI inequality which extends the HWI inequality established in
\cite{otto:generalization} and
\cite{CMV}, with the additional ``B'' referring here to the new barycentric
term. 
\begin{theorem} 
\label{theo4.1} {\rm (HWBI inequality)} 
Under the hypothesis of Theorem 2.1, 
we have for all $\rho_0, \rho_1\in {\cal P}_c(\Omega)$,
satisfying
$\mbox{supp}\,\rho_0\subset \Omega$,  and $P_F(\rho_0)\in
W^{1,\infty}(\Omega)$,  
\begin{equation}
\label{eqn:3.13}
{\rm H}^{F, W}_U(\rho_0|\rho_1)\leq
W_2(\rho_0,\rho_1)\sqrt{I_2(\rho_0|\rho_U)}
-\frac{\mu+\nu}{2}W_2^2(\rho_0,\rho_1)+\frac{\nu}{2}|{\rm b}(\rho_0)-{\rm
b}(\rho_1)|^2.
\end{equation}
 \end{theorem}
 \noindent{\bf Proof:}
Rewrite
(\ref{eqn:3.16}) as
\begin{equation}
\label{eqn:3.19}
 {\rm
H}_U^{F,W}(\rho_0|\rho_1)+\frac{\mu+\nu}{2}W_2^2(\rho_0,\rho_1) -
\frac{\nu}{2}|{\rm b}(\rho_0)-{\rm b}(\rho_1)|^2 
 \leq \frac{1}{2\sigma} W_2^2(\rho_0,\rho_1) + \frac{\sigma}{2}
I_2(\rho_0|\rho_U).
\end{equation}
Now minimize the right hand side of (\ref{eqn:3.19}) over $\sigma>0$.
The minimum is obviously achieved at $\bar \sigma =
\frac{W_2(\rho_0,\rho_1)}{\sqrt{I_2(\rho_0|\rho_U)}}$. This yields
(\ref{eqn:3.13}).\\

\noindent Setting $W=0$ (and then $\nu=0$) in Theorem \ref{theo4.1}, we obtain in
particular, the following HWI inequality first established by Otto-Villani
\cite{otto:generalization} in the case of the classical entropy $F(x)=x\ln
x$, and extended later on, for generalized entropy functions $F$ by
Carrillo, McCann and Villani in \cite{CMV}.

\begin{corollary}
\label{coro4.1}
{\rm (HWI inequalities \cite{CMV})}  
Under the hypothesis on $\Omega$ and $F$ in Theorem \ref{theo2.1}, let
$U:\R^n\rightarrow\R$ be a $C^2$-function with $D^2U\geq \mu I$, where
$\mu\in \R$. Then we have for all  $\rho_0, \rho_1\in {\cal P}_c(\Omega)$
satisfying $\mbox{supp}\,\rho_0\subset\Omega$,  and
$P_F(\rho_0)\in W^{1,\infty}(\Omega)$,
\begin{equation}
\label{eqn:3.20}
{\rm H}^F_{U}(\rho_0|\rho_1)\leq W_2(\rho_0,\rho_1)\sqrt{I(\rho_0|\rho_U)}
- \frac{\mu}{2}W_2^2(\rho_0,\rho_1).
\end{equation}
  \end{corollary}
If $U+W$ is uniformly convex (i.e., $\mu+\nu>0$) inequality
(\ref{eqn:3.16})  yields the
following extensions of the Log-Sobolev inequality:
\begin{corollary}
\label{coro4.2}
{(\rm Log-Sobolev inequalities with interaction potentials)}\\
 In addition to the hypothesis on $\Omega$, $F$, $U$ and $W$ in Theorem
\ref{theo2.1},
 assume $\mu+\nu>0$. Then for all  $\rho_0, \rho_1\in {\cal P}_c(\Omega)$
satisfying $\mbox{supp}\,\rho_0\subset\Omega$,  and $P_F(\rho_0)\in
W^{1,\infty}(\Omega)$,
we have
\begin{equation}
\label{eqn:3.25}
{\rm H}_U^{F,W}(\rho_0|\rho_1) - \frac{\nu}{2}|{\rm b}(\rho_0)-{\rm
b}(\rho_1)|^2
\leq\frac{1}{2(\mu+\nu)} I_2(\rho_0|\rho_U).
\end{equation}
In particular, if $b(\rho_0)=b(\rho_1)$, we have that
\begin{equation}
\label{eqn:3.26}
{\rm H}^{F,W}_U(\rho_0|\rho_1)\leq \frac{1}{2(\mu+\nu)}I_2(\rho_0|\rho_U).
\end{equation}
Furthermore, if $W$ is convex, then we have the following inequality,
established in \cite{CMV}
\begin{equation}
\label{eqn:3.27}
{\rm H}_U^{F,W}(\rho_0|\rho_1) \leq \frac{1}{2\mu} I_2(\rho_0|\rho_U).
\end{equation}
\end{corollary}
\noindent{\bf Proof:} (\ref{eqn:3.25}) follows easily from
(\ref{eqn:3.16}) by choosing $\sigma=\frac{1}{\mu+\nu}$, and (\ref{eqn:3.27}) follows
 from (\ref{eqn:3.25}), using $\nu=0$ because $W$ is convex.\\
In particular, setting $W=0$ in Corollary \ref{coro4.2}, one obtains the
following
generalized Log-Sobolev inequality obtained in \cite{CJMTU}, and in
\cite{cordero:inequalities} for generalized cost functions.
\begin{corollary}
\label{coro4.3}
{\rm (Generalized Log-Sobolev inequalities \cite{CJMTU},
\cite{cordero:inequalities})}\\
Assume that $\Omega$ and $F$ satisfy the assumptions in Theorem
\ref{theo2.1}, and that $U:\R^n\rightarrow \R$ is a $C^2$- uniformly
convex function with $D^2U\geq \mu I$, where $\mu>0$. Then for all
 $\rho_0, \rho_1\in {\cal P}_c(\Omega)$ satisfying
$\mbox{supp}\,\rho_0\subset\Omega$,  and $P_F(\rho_0)\in
W^{1,\infty}(\Omega)$, we have
\begin{equation}
\label{eqn:3.28}
{\rm H}_U^F(\rho_0|\rho_1) \leq\frac{1}{2\mu} I_2(\rho_0|\rho_U).
\end{equation}
 \end{corollary}
One can also deduce the following generalization of Talagrand's
inequality. We  note in
particular that when $W=0$, the result below is obtained previously by
Blower
\cite{blower}, Otto-Villani \cite{otto:generalization} and Bobkov-Ledoux
\cite{bobkov} for
the Tsallis entropy $F(x)=x\ln x$, and by Carrillo-McCann-Villani
\cite{CMV} for generalized
entropy functions $F$.
\begin{corollary}
\label{coro:4.4} {(\rm Generalized Talagrand Inequality with interaction
potentials)}\\
 In addition to the hypothesis on $\Omega$, $F$, $U$ and $W$ in Theorem
\ref{theo2.1},
assume $\mu+\nu>0$.  Then for all probability densities  $\rho$ on
$\Omega$,  we have
\begin{equation}
\label{eqn:3.30}
\frac{\nu+\mu}{2} W^2_2(\rho,\rho_U)-\frac{\nu}{2}|{\rm b}(\rho)-{\rm
b}(\rho_U)|^2
\leq  {\rm H}_U^{F,W}(\rho|\rho_U).
\end{equation}
 In particular, if $b(\rho)=b(\rho_U)$,
 we have that
\begin{equation}
\label{eqn:3.31}
W_2(\rho,\rho_U)\leq \sqrt{\frac{2{\rm H}^{F,W}_U(\rho|\rho_U)}{\mu+\nu}}.
\end{equation}
Furthermore, if $W$ is convex, then the following inequality established
in \cite{CMV}
holds:
\begin{equation}
\label{eqn:3.32}
W_2(\rho,\rho_U)\leq \sqrt{\frac{2{\rm H}^{F,W}_U(\rho|\rho_U)}{\mu}}.
\end{equation}
\end{corollary}
\noindent{\bf Proof:} (\ref{eqn:3.30}) follows from (\ref{eqn:3.16}) if we
use $\rho_0:=\rho_U$, $\rho_1:=\rho$, notice that $I_2(\rho_U|\rho_U)=0$, and
then let $\sigma$ go to $\infty$. (\ref{eqn:3.32}) follows from (\ref{eqn:3.30}), 
where we use $\nu=0$ because $W$ is convex.

\subsection{Inequalities with Boltzmann reference measures}

To each confinement potential $U:\R^n\rightarrow\R$ with $D^2U\geq
\mu I$ where $\mu\in \R$, one associates a Boltzmann reference measure
denoted by $\rho_U$ which is the normalized   $\frac{e^{-U}}{\sigma_U}$, where
$\sigma_U=\int_{\R^n}e^{-U}\,\mbox{d}x$ is assumed to be finite. To deduce
inequalities involving such reference measures, we can apply Proposition
\ref{prop4.1} with 
$F(x)=x\ln x$ and $W=0$ to get Gross' Log-Sobolev inequality (when
$U(x)=\frac{1}{2}|x|^2$) and its extension by Bakry and Emery in
\cite{bakry} (when $U$ is uniformly convex). We first state the following
HWI-type inequality from which we deduce
Otto-Villani's HWI inequality \cite{otto:generalization}, and the
Log-Sobolev
inequality of Gross \cite{Gross} and Bakry-Emery \cite{bakry}.

\begin{corollary}
\label{coro4.5}
Let $U:\R^n\rightarrow\R$ be a $C^2$-function with $D^2U\geq \mu I$ where
$\mu\in \R$.
Then for any
$\sigma>0$, the following holds for any nonnegative function
$f$ such that $f\rho_U\in W^{1,\infty}(\R^n)$ and
$\int_{\R^n}f\rho_U\,\mbox{d}x=1$:
\begin{equation}
\label{eqn:3.33}
\int_{\R^n}f\ln(f)\,\rho_U \mbox{d}x+\frac{1}{2}(\mu-\frac{1}{\sigma})
W_2^2(f\rho_U,\rho_U)\leq
\frac{\sigma}{2}
\int_{\R^n}\frac{|\,\nabla f\,|^2}{f}\,\rho_U \mbox{d}x.
\end{equation}
\end{corollary}
\noindent{\bf Proof:} First assume that $f$ has compact support, and set
$F(x)=x\ln x$, $\rho_0=f\rho_U,\,\rho_1=\rho_U$ and $W=0$ in
(\ref{eqn:3.16}). We
have that
\begin{equation}
\label{eqn:3.34}
H^F_U(f\rho_U|\rho_U)+\frac{1}{2}(\mu-\frac{1}{\sigma})W^2_2(f\rho_U,\rho_U)
\leq\frac{\sigma}{2}\int_{\R^n}\Big|\,\frac{\nabla(f\rho_U)}{f\rho_U}+U\,\Big|^2f\rho_U\,\mbox{d}x.
\end{equation}
By direct computations,
\begin{equation}
\label{eqn:3.35}
\frac{\nabla(f\rho_U)}{f\rho_U}=\frac{\nabla f}{f}-\nabla U,
\end{equation}
and
\begin{eqnarray}
\label{eqn:3.36}
{\rm H}^{F,W}_U(f\rho_U|\rho_U) &\leq&
\int_{\R^n}\left[f\rho_U\ln(f\rho_U)+Uf\rho_U-\rho_U\ln\rho_U-U\rho_U\right]\,\mbox{d}x\\
&=& \int_{\R^n}(f\rho_U\ln
f)\,\mbox{d}x+\ln\sigma_U\int_{\R^n}\left(\rho_U-f\rho_U\right)\,\mbox{d}x\nonumber\\
&=& \int_{\R^n}f\ln(f)\rho_U\,\mbox{d}x.\nonumber
\end{eqnarray}
Combining (\ref{eqn:3.34}) - (\ref{eqn:3.36}), we get (\ref{eqn:3.33}). We
finish the proof using a standard approximation argument.

\begin{coro}
\label{coro4.5'}
{(\rm Otto-Villani's HWI inequality \cite{otto:generalization})} 
Let $U:\R^n\rightarrow\R$ be a $C^2$-uniformly convex function with
$D^2U\geq \mu I$, where $\mu>0$. Then, for any nonnegative
function $f$ such that $f\rho_U\in W^{1,\infty}(\R^n)$ and
$\int_{\R^n}f\rho_U\,\mbox{d}x=1$,
\begin{equation}
\label{eqn:3.36'}
\int_{\R^n} f\ln(f)\rho_U\,\mbox{d}x \leq
W_2(\rho_U,f\rho_U)\sqrt{I(f\rho_U|\rho_U)}-\frac{\mu}{2}W_2^2(f\rho_U,\rho_U),
\end{equation} 
where
\[ I(f\rho_U|\rho_U)=\int_{\R^n}\frac{|\,\nabla
f\,|^2}{f}\rho_U\,\mbox{d}x.\]
\end{coro}
\noindent{\bf Proof:} It is similar to the proof of Theorem \ref{theo4.1}.
Rewrite
(\ref{eqn:3.33}) as
\[ \int_{\R^n} f\ln(f)\rho_U\,\mbox{d}x +\frac{\mu}{2}
W^2_2(f\rho_U,\rho_U)\leq \frac{\mu}{2\sigma}
W^2_2(f\rho_U,\rho_U)+\frac{\sigma}{2} I(f\rho_U|\rho_U),\]
and show that the minimum over $\sigma>0$ of the right hand side is
attained at
$\bar{\sigma}=\frac{W_2(f\rho_U,\rho_U)}{\sqrt{I(f\rho_U|\rho_U)}}$.\\
\noindent Setting $f:=g^2$ and $\sigma:=\frac{1}{\mu}$ in (\ref{eqn:3.36'}),
one obtains
the following extension of Gross' \cite{Gross} Log-Sobolev inequality
first
established by Bakry and Emery in \cite{bakry}.

\begin{corollary} {\rm (Original Log Sobolev inequality \cite{bakry},
\cite{Gross})} 
\label{coro4.6}
Let $U:\R^n\rightarrow \R$ be a $C^2$-uniformly convex function with
$D^2U\geq \mu I$ where $\mu>0$.
Then, for
any function $g$ such that $g^2\rho_U\in W^{1,\infty}(\R^n)$ and
$\int_{\R^n}g^2\rho_U\,\mbox{d}x=1$, we have
\begin{equation}
\label{eqn:3.37}
\int_{\R^n}g^2\ln(g^2)\,\rho_U\mbox{d}x \leq\frac{2}{\mu}\int_{\R^n}
|\,\nabla g\,|^2 \,\rho_U \mbox{d}x.
\end{equation}
\end{corollary}
As pointed out by Rothaus in \cite{Rothaus}, the above Log-Sobolev
inequality implies the Poincar\'e's inequality.
\begin{coro}
\label{coro4.6.1}
{(\rm Poincar\'e's inequality)}
Let $U:\R^n\rightarrow \R$ be a $C^2$-uniformly convex function with
$D^2U\geq \mu I$ where $\mu>0$. Then, for any function $f$ such
that $f\rho_U\in W^{1,\infty}(\R^n)$ and
$\int_{\R^n}f\rho_U\,\mbox{d}x=0$, we have
\begin{equation}
\label{eqn:p1}
\int_{\R^n}f^2\rho_U\,\mbox{d}x \leq \frac{1}{\mu}\int_{\R^n}|\,\nabla
f\,|^2\rho_U\,\mbox{d}x.
\end{equation}
\end{coro}
\noindent{\bf Proof:} From (\ref{eqn:3.37}), we have that
\begin{equation}
\label{eqn:p2}
\int_{\R^n}f_\epsilon\ln(f_\epsilon)\,\rho_U\,\mbox{d}x \leq
\frac{1}{2\mu}\int_{\R^n}\frac{|\,\nabla
f_\epsilon\,|^2}{f_\epsilon}\rho_U\,\mbox{d}x,
\end{equation}
where $f_\epsilon=1+\epsilon f$ for some $\epsilon>0$. Using that
$\int_{\R^n}f\rho_U\,\mbox{d}x=0$, we have for small $\epsilon$,
\begin{equation}
\label{eqn:p3}
\int_{R^n} f_\epsilon\ln(f_\epsilon)\rho_U\,\mbox{d}x =
\frac{\epsilon^2}{2} \int_{\R^n} f^2\rho_U\,\mbox{d}x + o(\epsilon^3),
\end{equation}
and
\begin{equation}
\label{eqn:p4}
\int_{\R^n}\frac{|\,\nabla f_\epsilon\,|^2}{f_\epsilon}\,\rho_U\,\mbox{d}x
= \epsilon^2\int_{\R^n}|\,\nabla f\,|^2\rho_U\,\mbox{d}x + o(\epsilon^3).
\end{equation}
We combine (\ref{eqn:p2}) - (\ref{eqn:p4}) to have that
\begin{equation}
\label{eqn:p5}
\int_{\R^n}f^2\rho_U\,\mbox{d}x \leq \frac{1}{\mu}\int_{\R^n} |\,\nabla
f\,|^2\rho_U\,\mbox{d}x + o(\epsilon).
\end{equation}
We let $\epsilon$ go to $0$ in (\ref{eqn:p5}) to conclude
(\ref{eqn:p1}).\\

\noindent If we apply Corollary \ref{coro:4.4} to $F(x)=x\ln x$ when $W=0$, we
obtain the
following extension of Talagrand's inequality established by Otto and
Villani in
\cite{otto:generalization}.

\begin{coro}
\label{coro4.6'}{\rm (Original Talagrand's inequality \cite{T},
\cite{otto:generalization})}
Let $U:\R^n\rightarrow \R$ be a $C^2$-uniformly convex function with
$D^2U\geq \mu I$ where $\mu>0$. Then, for any nonnegative
function $f$ such that $\int_{\R^n}f\rho_U\,\mbox{d}x=1$, we have
\begin{equation}
\label{eqn:3.37.1}
W_2(f\rho_U,\rho_U)\leq
\sqrt{\frac{2}{\mu}\int_{\R^n}f\ln(f)\rho_U\,\mbox{d}x}.
\end{equation}
\end{coro}
In particular, if $f=\frac{\1_B}{\rho_U(B)}$ for some measurable subset
$B$ of $\R^n$, where $d\gamma(x)=\rho_U(x)\mbox{d}x$ and $\1_B$ is the
characteristic function of $B$, one obtains the following inequality in the
concentration of measures in Gauss space, first proved by Talagrand
building on an argument by Marton (see details in Villani
\cite{villani:topics}).  

\begin{corollary} {\rm (Concentration of measure inequality)} 
\label{coro4.6''}
Let $U:\R^n\rightarrow \R$ be a $C^2$-uniformly convex function with
$D^2U\geq \mu I$ where $\mu>0$. Then, for any
$\epsilon$-neighborhood $B_\epsilon$ of a measurable set $B$ in $\R^n$, we
have
 \begin{equation}
\label{eqn:3.37.2}
  \gamma (B_\epsilon) \geq 1- e^{-\frac{\mu}{2}\left(\epsilon-\sqrt
{\frac{2}{\mu}\ln\left(\frac{1}{\gamma (B)}\right)}\right)^2},
\end{equation}
where $\epsilon\geq
\sqrt{\frac{2}{\mu}\ln\left(\frac{1}{\gamma(B)}\right)}$.
\end{corollary}
\noindent{\bf Proof:} Using $f=f_B=\frac{\1_B}{\gamma(B)}$ in
(\ref{eqn:3.37.1}), we have that
\[ W_2(f_B\rho_U,\rho_U)\leq
\sqrt{\frac{2}{\mu}\ln\left(\frac{1}{\gamma(B)}\right)},\]
and then, we obtain from the triangle inequality that
\begin{equation}
\label{eqn:3.37.3}
 W_2(f_B\rho_U,f_{\R^n\backslash B_\epsilon}\rho_U)\leq
\sqrt{\frac{2}{\mu}\ln\left(\frac{1}{\gamma(B)}\right)}
+\sqrt{\frac{2}{\mu}\ln\left(\frac{1}{1-\gamma(B_\epsilon)}\right)}.
\end{equation}
But since $|\,x-y\,|\geq \epsilon$ for all $(x,y)\in B\times
(\R^n\backslash B_\epsilon)$, we have that
\begin{equation}
\label{eqn:3.37.4}
W_2(f_B\rho_U,\rho_U) \geq \epsilon.
\end{equation}
We combine (\ref{eqn:3.37.3}) and (\ref{eqn:3.37.4}) to deduce that
\[\ln\left(\frac{1}{1-\gamma(\R^n\backslash B_\epsilon)}\right) \geq
\frac{\mu}{2}\left(\epsilon-\sqrt{\frac{2}{\mu}\ln\left(\frac{1}{\gamma(B)}\right)}\right)^2,\]
which leads to (\ref{eqn:3.37.2}).

\section{Trends to equilibrium}
\label{sect5}

We use Corollary \ref{coro4.3} and Corollary \ref{coro:4.4} to
recover rates
of convergence for solutions to equation

\begin{equation}
\label{eqn:5.2}
\left\{\begin{array}{lcl}
\frac{\partial\rho}{\partial t} = \mbox{div}\left\{\rho\nabla
\left(F'(\rho)+V+W\star\rho\right)\right\} &\mbox{in}&
(0,\infty)\times\R^n\\ \\
\rho(t=0)=\rho_0 &\mbox{in}& \{0\}\times\R^n,
\end{array}\right.
\end{equation}
recently shown by Carrillo, McCann and Villani in \cite{CMV}. Here we
consider
 the case where $V+W$ is uniformly convex and $W$ convex, and the case
when only
$V+W$ is uniformly convex but the barycentre $b\left(\rho(t)\right)$ of
any solution
$\rho(t,x)$ of (\ref{eqn:5.2}) is invariant in $t$. For a background and
other
cases of convergence to equilibrium for this equation, we refer to
\cite{CMV} and
the references therein.

\begin{coro}
\label{prop5.3}
{(\rm Trend to equilibrium)}
Let $F:[0,\infty)\rightarrow\R$ be strictly convex, differentiable on
$(0,\infty)$ and satisfies $F(0)=0$,
$\lim_{x\rightarrow\infty}\frac{F(x)}{x}=\infty$, and
$x\mapsto x^nF(x^{-n})$ is convex and non-increasing. Let
$V,\,W:\R^n\rightarrow
[0,\infty)$ be  respectively $C^2$-confinement and interaction potentials
with $D^2V\geq
\lambda I$ and $D^2W\geq \nu I$, where $\lambda, \nu\in\R$. Assume that
the initial
probability density $\rho_0$ has finite total energy. Then
\begin{enumerate}
\item If $V+W$ is uniformly convex (i.e., $\lambda+\nu>0$) and $W$ is
convex 
(i.e. $\nu\geq 0$), then, for any solution $\rho$ of (\ref{eqn:5.2}), such
that ${\rm
H}^{F,W}_V\left(\rho(t)\right)<\infty$, we have:
\begin{equation}
\label{eqn:5.13}
{\rm H}^{F,W}_V\left(\rho(t)|\rho_V\right)\leq e^{-2\lambda t}{\rm
H}^{F,W}_V(\rho_0|\rho_V),
\end{equation}
and
\begin{equation}
\label{eqn:5.14}
W_2\left(\rho(t),\rho_V\right)\leq e^{-\lambda t}\sqrt{\frac{2{\rm
H}^{F,W}_V(\rho_0|
\rho_V)}{\lambda}}.
\end{equation}
 \item If $V+W$ is uniformly convex (i.e., $\lambda+\nu> 0$) and if we
assume
that the barycentre $b\left(\rho(t)\right)$ of any solution $\rho(t,x)$ of
(\ref{eqn:5.2})
is invariant in $t$, then, for any solution $\rho$ of (\ref{eqn:5.2}) such
that ${\rm
H}^{F,W}_V\left(\rho(t)\right)<\infty$, we have:
\begin{equation}
\label{eqn:5.15}
{\rm H}^{F,W}_V\left(\rho(t)|\rho_V\right)\leq e^{-2(\lambda+\nu)t}{\rm
H}^{F,W}_V(\rho_0|\rho_V),
\end{equation}
and
\begin{equation}
\label{eqn:5.16}
W_2\left(\rho(t),\rho_V\right)\leq e^{-2(\lambda+\nu)t}\sqrt{\frac{2{\rm
H}^{F,W}_V(\rho_0|\rho_V)}{\lambda+\nu}}.
\end{equation}
\end{enumerate}
\end{coro}
\noindent{\bf Proof:} Under the assumptions on $F$, $V$ and $W$ in
Corollary
\ref{prop5.3}, it is known (see \cite{CMV}, and references therein) that
the total
energy ${\rm H}^{F,W}_V$ -- which is a Lyapunov functional associated with
(\ref{eqn:5.2}) -- has a unique minimizer $\rho_V$ defined by
\[
      \rho_V\nabla\left(F^\prime(\rho_{_V})+V+W\star \rho_{_V}\right) =
0\quad \mbox{a.e.}
\]
Now,
let $\rho$ be a -- smooth -- solution of (\ref{eqn:5.2}). We have the
following energy
dissipation equation
\begin{equation}
\label{eqn:5.17}
\frac{d}{dt}\,{\rm
H}^{F,W}_V\left(\rho(t)|\rho_V\right)=-I_2\left(\rho(t)|\rho_V\right).
\end{equation}
Combining (\ref{eqn:5.17}) with (\ref{eqn:3.27}), we have that
\begin{equation}
\label{eqn:5.18}
\frac{d}{dt}\,{\rm H}^{F,W}_V\left(\rho(t)|\rho_V\right) \leq -2\lambda
{\rm H}^{F,W}_V\left(\rho(t)|\rho_V\right).
\end{equation}
We integrate (\ref{eqn:5.18}) over $[0,t]$ to conclude (\ref{eqn:5.13}).
(\ref{eqn:5.14}) follows directly from (\ref{eqn:3.32}) and
(\ref{eqn:5.13}).\\
To prove (\ref{eqn:5.15}), we use (\ref{eqn:5.17}) and (\ref{eqn:3.26}) to
have that\begin{equation}
\label{eqn:5.19}
\frac{d}{dt}{\rm H}^{F,W}_V\left(\rho(t)|\rho_V\right)\leq
-2(\lambda+\nu){\rm H}^{F,W}_V\left(\rho(t)|\rho_V\right).
\end{equation}
We integrate (\ref{eqn:5.19}) over $[0,t]$ to conclude (\ref{eqn:5.15}).
As before, (\ref{eqn:5.16}) is a consequence of (\ref{eqn:5.15}) and
(\ref{eqn:3.31}).\\

\noindent We now apply Corollary \ref{prop5.3} to obtain rates of convergence
to equilibrium for some  equations of the form (\ref{eqn:5.2}) studied in the
literature by many authors. 

\begin{trivlist}
\item $\bullet$ \ If $W=0$ and $F(x)=x\ln x$ in which case (\ref{eqn:5.2}) is
the linear
Fokker-Planck equation $\frac{\partial\rho}{\partial t}=\Delta\rho
+\mbox{div}(\rho\nabla
V)$, Corollary \ref{prop5.3} gives an exponential decay in relative
entropy of solutions of
this equation to the Gaussian density
$\rho_V=\frac{e^{-V}}{\sigma_V},\;\sigma_V=\int_{\R^n}e^{-V}\,\mbox{d}x$,
at the rate
$2\lambda$ when $D^2V\geq \lambda I$ for some $\lambda>0$, and an
exponential decay in the
Wasserstein distance, at the rate $\lambda$.
\item $\bullet$ \ If $W=0$, $F(x)=\frac{x^m}{m-1}$ where $1\neq m\geq
1-\frac{1}{n}$, and
$V(x)=\lambda\frac{|\,x\,|^2}{2}$ for some $\lambda>0$, in which case
(\ref{eqn:5.2}) is the
rescaled porous medium equation ($m>1$), or fast diffusion equation
$(1-\frac{1}{n}\leq
m<1)$, that is  $\frac{\partial\rho}{\partial t} =\Delta\rho^m
+\mbox{div}(\lambda x\rho)$,
Corollary \ref{prop5.3} gives an exponential decay in relative entropy of
solutions of this
equation to the Barenblatt-Prattle profile
$\rho_V(x)=\left[\left(C+\frac{\lambda(1-m)}{2m}|\,x\,|^2\right)^{\frac{1}{m-1}}\right]^+$
(where $C>0$ is such that $\int_{\R^n}\rho(x)\,\mbox{d}x=1$) at the rate
$2\lambda$, and an
exponential decay in the Wasserstein distance at the rate $\lambda$.
\end{trivlist}

\section{The Energy-Entropy production Duality Formula}
\label{sect6}

In this section, we apply Theorem \ref{theo2.1} with $V=W=0$, to obtain the
following intriguing duality formula.

\begin{proposition} {\rm (The Energy-Entropy Duality Formula)} Under the hypothesis
of Theorem 2.1, we have for any
$\rho_0,
\rho_1 \in {\cal P}_c(\Omega)$  satisfying
$\mbox{supp}\,\rho_0\subset \Omega$ and $P_F(\rho_0)\in
W^{1,\infty}(\Omega)$, and any Young function $c:\R^n\to\R$:
\begin{equation}
-{\rm H}^F_{c}(\rho_1)\leq -{\rm H}^{F+nP_F}(\rho_0)+\int_\Omega \rho_0
c^\star\left(-\nabla(F^\prime\circ\rho_0)\right)\,\mbox{d}x.
\end{equation}
 Moreover, equality holds whenever $\rho_0=\rho_1=\rho_c$ where $\rho_c$ is
a probability
density on $\Omega$ such that $\nabla(F'(\rho_c)+c)=0$ a.e.
\end{proposition} 
Motivated by the recent work of Cordero-Nazaret-Villani \cite{cordero:mass}, 
we show that this inequality points to a remarkable correspondence between 
ground state solutions of some quasilinear PDEs or semi-linear equations which
appear as Euler-Lagrange equations of the entropy production functionals, and
stationary solutions of Fokker-Planck type equations. 
 
\begin{corollary}
\label{coro3.1}
Under the hypothesis of Theorem 2.1, let $\psi:\R\rightarrow [0,\infty)$
differentiable be chosen in such a
way that $\psi(0)=0$ and $|\,\psi^{\frac{1}{p}} (F' \circ \psi)'\,| = K$
where $p>1$, and
$K$ is chosen to be 1 for simplicity. Then, for any Young function $c$
with $p$-homogeneous
Legendre transform $c^*$, we have the following inequality:
     \begin{equation}
\label{eqn:3.1}
        \sup\{ - \int_{\Omega} F(\rho) + c\rho; \rho  \in {\cal P}_c
(\Omega)\}
\leq
\inf \{ \int_{\Omega} c^{*} (-\nabla f) - G_F \circ \psi(f); f\in
C^\infty_0(\Omega),
\int_\Omega \psi (f) =1\}
     \end{equation}
where $G_F(x) := (1-n)F(x) +nxF'(x)$.\\
 Furthermore,  equality holds in
(\ref{eqn:3.1}) if
    there exists $\bar f$ (and $\bar\rho=\psi(\bar f)$) that satisfies
     \begin{equation}
\label{eqn:3.2}
      -(F' \circ \psi)'(\bar f) \nabla \bar f (x) = \nabla  c(x)\;\; a.e.
\end{equation}
Moreover, $\bar f$ solves
       \begin{equation}
\label{eqn:3.3}
\begin{array}{ll}
       \mbox{\rm div} \{ \nabla c^{*} (-\nabla f) \} -
        (G_F \circ \psi)'(f) =
       \lambda \psi'(f) & \mbox{\rm in $\Omega$}\\
                   \nabla c^{*} (-\nabla f)
                   \cdot \nu = 0 & \mbox{\rm on $\partial\Omega$},
                   \end{array}
          \end{equation}
for some $\lambda\in\R$, while $\bar{\rho}$ is a stationary solution of

\begin{equation}
\label{eqn:3.3'}
\begin{array}{ll}
       \frac{\partial \rho}{\partial t}=\mbox{\rm div}
\{\rho\nabla\left(F'(\rho)+c\right) \} & \mbox{\rm in
$(0,\infty)\times\Omega$}\\
               \rho\nabla\left(F'(\rho)+c\right)\cdot \nu = 0 & \mbox{\rm
on $(0,\infty)\times\partial\Omega$}.
                   \end{array}
\end{equation}

\end{corollary}
\noindent {\bf Proof:}
Assume that $c^*$ is $p$-homogeneous, and let $Q''(x) =
x^{\frac{1}{q}} F''(x)$ where $q$ is the conjugate of $p$.
     Let \[
          J(\rho):= - \int_{\Omega} [ F(\rho(y)) + c(y) \rho(y) ] dy
       \]
       and
       \[
          {\tilde J}(\rho):= -\int_{\Omega} (F+ n P_F) (\rho(x)) dx
                  + \int_{\Omega} c^{*} ( -\nabla ( Q'(\rho(x))) dx.
       \]
Equation (\ref{eqn:2.1}) (where we use $V=W=0$, and then $\lambda=\nu=0$)
then becomes
        \begin{equation}
        \label{eqn:3.4}
          J(\rho_1) \leq {\tilde J}(\rho_0)
       \end{equation}
       for all probability densities $\rho_0, \rho_1$ on $\Omega$ such
that
$\mbox{supp}\,\rho_0\subset\Omega$ and $P_F(\rho_0)\in
W^{1,\infty}(\Omega)$. If
$\bar \rho$ satisfies
       \[
          -\nabla (F'(\bar \rho(x) ))= \nabla \mbox{$c(x)$ a.e.},
       \]
then equality holds in (\ref{eqn:3.4}), and $\bar \rho$ is an extremal
of the variational
    problems
     \[
          \sup \{ J(\rho);\  \rho  \in {\cal P}_c (\Omega)\}=
          \inf \{ {\tilde J}(\rho); \rho  \in {\cal P}_c (\Omega),
\mbox{supp}\,\rho\subset\Omega,  P_F(\rho)\in W^{1,\infty}(\Omega)  \}.
           \]
     In particular, $\bar\rho$ is a solution of
\begin{equation}
\label{eqn:3.5}
\begin{array}{ll}
                   \mbox{\rm div} \{ \rho \nabla (F'(\rho)+c) \} = 0 &
                   \mbox{\rm in $\Omega$}\\
                   \rho \nabla (F'(\rho)+c) \cdot \nu = 0 & \mbox{\rm on
                   $\partial \Omega$}.
                   \end{array}
          \end{equation}
Suppose now $\psi:\R \rightarrow [0,\infty)$ differentiable, $\psi
(0)=0$ and that $\bar f\in C_0^\infty(\Omega)$ satisfies
$-(F' \circ \psi)'(\bar f) \nabla \bar f (x) = \nabla \mbox{$c(x)$
a.e.}$ Then equality holds in (\ref{eqn:3.4}), and $\bar f$ and $\bar \rho
=\psi (\bar f)$ are
extremals of the
following variational problems
     \[
          \inf \{ I(f);\ f\in C^{\infty}_0(\Omega),  \int_\Omega\psi (f)=
1 \}=
          \sup \{ J(\rho); \rho  \in {\cal P}_c (\Omega) \}
           \]
where

\[
I (f) ={\tilde J} (\psi (f))=-\int_{\Omega} [ F \circ \psi + n P_F
\circ \psi ] (f)
                  + \int_{\Omega} c^{*} ( -\nabla ( Q' \circ \psi
                  (f) )).
\]
If now $\psi$ is such that $|\,\psi^{\frac{1}{p}} (F' \circ \psi)'\,|
= 1 $, then $|\,(Q' \circ \psi)'\,| = 1$ and
\[
I (f) = -\int_{\Omega} [ F \circ \psi + n P_F \circ \psi ] (f)
                  + \int_{\Omega} c^{*} ( -\nabla f )),
\]
because $c^*$ is $p$-homogeneous. This proves (\ref{eqn:3.1}). The
Euler-Lagrange equation of the variational problem
\[ \inf \Big\{ \int_{\Omega} c^{*} (-\nabla (f)) -[ F \circ \psi + n P_F
\circ
\psi ] (f) ;\ \int_\Omega \psi (f) =1\Big\}
       \]
reads as
\begin{equation}
\label{eqn:3.6}
     \begin{array}{ll}
       \mbox{div} \{ \nabla c^{*} (-\nabla f) \} -
        (G_F \circ \psi)'(f) =
       \lambda \psi'(f) & \mbox{\rm in $\Omega$}\\
                   \nabla c^{*} (-\nabla f)
                   \cdot \nu = 0 & \mbox{\rm on $\partial\Omega$}
                   \end{array}
          \end{equation}
where $\lambda\in\R$ is a Lagrange multiplier, and $G(x) = (1-n) F(x) + nx
F'(x)$. This proves (\ref{eqn:3.3}). To prove that the maximizer
$\bar{\rho}$ of
\[ \sup\{-\int_\Omega \left(F(\rho)+c\rho\right)\,\mbox{d}x;\; \rho\in
{\cal P}_c(\Omega)\}\]
is a stationary solution of (\ref{eqn:3.3'}), we refer to \cite{JKO} and
\cite{O:geometry}.\\
Now, we apply Corollary \ref{coro3.1} to the functions $F(x)=x\ln x,
\psi(x)=|\,x\,|^p$ and $c(x)=(p-1)|\,\mu x\,|^q$, with $\mu>0$ and
$c^*(x)=\frac{1}{p}\Big|\,\frac{x}{\mu}\,\Big|^p$ and
$\frac{1}{p}+\frac{1}{q}=1$, to derive a duality between stationary
solutions of Fokker-Planck equations, and ground state solutions of some
semi-linear equations.  We note here that the condition
$|\,\psi^{\frac{1}{p}}(F'\circ\psi)'\,|=K$ holds for $K=p$. We obtain the
following:
\begin{coro}
\label{coro3.2}
Let $p>1$ and let $q$ be its conjugate ($\frac{1}{p}+\frac{1}{q}=1$).
For all $f\in W^{1,p}(\R^n)$, such that $\|\,f\,\|_p =1$, any
probability density $\rho$ such that $\int_{\R^n}\rho(x) |x|^q
dx<\infty$, and any $\mu>0$, we have
\begin{equation}
\label{eqn:3.7}
J_\mu(\rho)\leq I_\mu(f),
\end{equation}
where
\[ J_\mu(\rho):=-\int_{\R^n}
\rho\ln\left(\rho\right)\,\mbox{d}y-(p-1)\int_{\R^n} |\,\mu
y\,|^q\rho(y)\,\mbox{d}y,\]
and
\[ I_\mu(f):=-\int_{\R^n} |\,f\,|^p\ln\left(|\,f\,|^p\right)
+\int_{\R^n}\Big|\,\frac{\nabla f}{\mu}\,\Big|^p -n.\]
Furthermore, if $h\in W^{1,p}(\R^n)$ is such that $h\geq 0$,
$\|\,h\,\|_p=1$, and
\[ \nabla h(x)=-\mu^{q} x |\,x\,|^{q-2} h(x)\quad \mbox{a.e.,}\]
then
\[ J_\mu(h^p)=I_\mu(h).\]
Therefore, $h$ (resp., $\rho=h^p$) is an extremum of the variational
problem:
\[
\qquad \sup\{\,J_\mu(\rho): \rho\in W^{1,1}(\R^n),\;
     \|\,\rho\,\|_1=1\}= \inf\{\,I_\mu(f): f\in W^{1,p}(\R^n),
     \|\,f\,\|_p=1\}.
\]
\end{coro}
It follows that $h$ satisfies the Euler-Lagrange equation corresponding to
the
constraint minimization problem, i.e., $h$ is a solution of
\begin{equation}
\label{eqn:3.8}
      \mu^{-p}\Delta_pf + p f|\,f\,|^{p-2}\ln(|\,f\,|)=\lambda
f|\,f\,|^{p-2},
\end{equation}
where $\lambda$ is a Lagrange multiplier. On the other hand, $\rho=h^p$ is
a stationary solution of the Fokker-Planck equation:
\begin{equation}
\label{eqn:3.8.1}
                 \frac{\partial u}{\partial t}= \Delta u +
\mbox{div}(p\mu^q|x|^{q-2}x u  ).
          \end{equation}
We can also apply Corollary \ref{coro3.1} to recover the duality
associated to the Gagliardo-Nirenberg inequalities obtained recently
in \cite{cordero:mass}.
\begin{coro}
\label{coro3.3}
Let $1<p<n$, and  $r\in\left(0,\frac{np}{n-p}\right]$ such that $r\neq p$.
Set
     $\gamma:=\frac{1}{r}+\frac{1}{q}$, where
$\frac{1}{p}+\frac{1}{q}=1$. Then, for $f \in W^{1,p}(\R^n)$ such
that $\|\,f\,\|_r=1$, for any probability density $\rho$ and for all
$\mu>0$, we have
\begin{equation}
\label{eqn:3.10}
J_\mu(\rho)\leq I_\mu(f)
\end{equation}
where
\[ J_\mu(\rho):=-\frac{1}{\gamma-1}\int_{\R^n}\rho^{\gamma}
-\frac{r\gamma\mu^{q}}{q}\int_{\R^n} |\,y\,|^{q}\rho(y)\,\mbox{d}y,\]
and
\[I_\mu(f):=-\left(\frac{1}{\gamma-1}+n\right) \int_{\R^n}
|\,f\,|^{r\gamma} +\frac{r\gamma}{p\mu^p} \int_{\R^n} |\,\nabla
f\,|^p.\]
Furthermore, if $h\in W^{1,p}(\R^n)$ is such that $h\geq 0$,
$\|\,h\,\|_r=1$, and
\[\nabla h(x)=-\mu^{q} x|\,x\,|^{q-2}h^{\frac{r}{p}}(x)\quad
\mbox{a.e.,}\]
     then
\[ J_\mu(h^r)=I_\mu(h).\]
Therefore, $h$ (resp., $\rho=h^r$) is an extremum of the variational
problems
\[
\qquad \sup\{\,J_\mu(\rho): \rho\in W^{1,1}(\R^n),\;
     \|\,\rho\,\|_1=1\}= \inf\{\,I_\mu(f): f\in W^{1,p}(\R^n),
     \|\,f\,\|_r=1\}.
\]
\end{coro}
\noindent{\bf Proof:} Again, the proof follows from Corollary
\ref{coro3.1}, by using now $\psi(x)=|\,x\,|^r$ and
$F(x)=\frac{x^\gamma}{\gamma-1}$, where $1\neq \gamma\geq 1-\frac{1}{n}$,
which follows from the fact that $p\neq r\in
\left(0,\frac{np}{n-p}\right]$. Indeed, for this value of $\gamma$,
the function $F$  satisfies the conditions of Corollary \ref{coro3.1}. The
Young function is now  $c(x)=\frac{r\gamma}{q}|\,\mu x\,|^q$, that is,
$c^*(x)=\frac{1}{p(r\gamma)^{p-1}}\Big|\,\frac{x}{\mu}\,\Big|^p$, and
the condition $|\,\psi^{\frac{1}{p}}(F'\circ\psi)'\,|=K$ holds with
$K=r\gamma$.\\
    Moreover, if $h\geq 0$ satisfies (\ref{eqn:3.2}), which is here,
\[ -\nabla h(x)=\mu^q x|\,x\,|^{q-2}h^{\frac{r}{p}}(x)\;\;\mbox{a.e.},\]
then $h$ is extremal in the minimization problem defined in Corollary
\ref{coro3.3}.\\
As above, we also note that  $h$ satisfies the Euler-Lagrange
equation corresponding to the constraint minimization problem, that
is, $h$ is a solution of
\begin{equation}
\label{eqn:3.11}
      \mu^{-p}\Delta_p f +
\left(\frac{1}{\gamma-1}+n\right)f|\,f\,|^{r\gamma-2}=\lambda
f|\,f\,|^{r-2},
\end{equation}
where $\lambda$ is a Lagrange multiplier. On the other hand,
$\rho=h^r$ is a stationary solution of the evolution equation:
\begin{equation}
\label{eqn:3.12}
     \frac{\partial u}{\partial t}= \Delta u^\gamma +
\mbox{div}(r\gamma\mu^q|x|^{q-2}x
u).
          \end{equation}

\noindent {\bf Example}: In particular, when $\mu=1, p=2,
\gamma=1-\frac{1}{n}$ and then $r=2^*=\frac{2n}{n-2}$ is the critical
Sobolev exponent, then Corollary \ref{coro3.3} yields a duality between
solutions of (\ref{eqn:3.11}), which is here the Yamabe equation:
\[-\Delta f=\lambda f|\,f\,|^{2^*-2},\]
(where $\lambda$ is the Lagrange multiplier due to the constraint
$\|\,f\,\|_{2^*}=1$), and stationary solutions of (\ref{eqn:3.12}), which
is here the rescaled fast diffusion equation:
\[\frac{\partial u}{\partial t}=\Delta
u^{1-\frac{1}{n}}+\mbox{div}\,\left(\frac{2n-2}{n-2}xu\right).\]

\end{document}